\documentclass[onecolumn,amsmath,showkeys,noshowpacs,amssymb,nofootinbib]{revtex4}

\usepackage{graphicx,color}

\usepackage{amsbsy}
\usepackage{amsthm,mathrsfs,amsopn}
\usepackage{multirow}

\begin{document}

\title{On the location and the strength of controllers to desynchronize coupled Kuramoto oscillators}
\author{Martin Moriam\'e$^{1}$ and Timoteo Carletti$^{1}$}

\affiliation{1. Department of Mathematics and Namur Institute for Complex Systems, naXys, University of Namur, rue Graf\'e 2, B 5000 Namur, Belgium}

\begin{abstract} 
	Synchronization is an ubiquitous phenomenon in dynamical systems of networked oscillators. While it is often a goal to achieve, in some context one would like to decrease it, e.g., although synchronization is essential to the good functioning of brain dynamics, hyper-synchronization can induce problems like epilepsy seizures. Motivated by this problem, scholars have developed pinning control schemes able to decrease synchronization in a system. Focusing on one of these methods, the goal of the present work is to analyse which is the best way to select the controlled nodes, i.e. the one that guarantees the lower synchronization rate. We show that hubs are generally the most advantageous nodes to control, especially when the degree distribution is heterogeneous. Nevertheless, pinning a too large number of hubs is in general not an appropriate choice. Our results are in line with previous works that studied pinning control aimed to increase synchronization. These observations shed light on an interesting universality of good practice of node selection disregarding the actual goal of the control scheme. 
\end{abstract}



\maketitle

\vspace{0.8cm}

\section{Introduction}
\label{sec:intro}

Synchronization is a widespread phenomenon in Nature, in particular in the living kingdom where interacting subsystems, being them biochemical cycles, cells or organs, are capable to synchronize their behaviors to eventually exhibit unison rhythms~\cite{pikovsky2003}. This behavior allows to increase the signal vs. noise ratio and also to positively react to changes in the external environment, ultimately determining suitable conditions to live~\cite{goldbeter_berridge_1996}.

On the other hand, living beings should be able to deal with several synchronized processes and thus, by suitably controlling them, they should be able to switch them on and off on demand. There are however cases where living being are unable to perform correctly such control task and hyper-synchronization occurs; this behavior can result in a pathological state, e.g., tremors in neurodegenerative diseases such as Parkinson’s disease (PD), or epileptic seizures. To overcome such issue most patients suffering from PD or epilepsy resort to long term drug treatments, whose results are only partially satisfying~\cite{Parkinson,Parkes:1974wr,Holloway:2005wh,Deckers:2000tz} and are often prone to  side effects. An alternative to chemical treatments is provided by the neurostimulation, that ultimately relies on the possibility to directly act on (bunch of) neurones by injecting a signal, i.e., a current capable to desynchronize their dynamics~\cite{neurostim,tass,neurocontrol}, among those methods the Deep Brain Stimulation (DBS) is among the most frequently used~\cite{DBS,Bronstein:2011hh,Theodore:2004gv}.

Scholars have thus devised suitable control techniques to reduce or even suppress hyper-synchronization in coupled oscillators. A particularly active research field is the one related to brain modeling because of the many possible medical applications. In this case however, one has to consider additional constraints, e.g., limiting the invasiveness of the method while still achieving suitable level of desynchronization.

Working with the abstract model of neurones given by the Kuramoto nonlinear oscillator~\cite{kuramoto1975,kuramoto}, it has been recently proposed a control scheme~\cite{gjata,AsllaniCarletti2018}, rooted on the Hamiltonian control theory~\cite{vittot,chandre,ciraolo,carletti}, that is minimally invasive in the number of required controllers but also in the strength of the injected signal. Indeed it has been compared with the Proportional-Differential feedback method~\cite{pyragas} and the results show that the former scheme exhibits the same efficiency as the latter, with the same number of controllers but with a weaker injected signal.

The efficiency of the above mentioned control scheme has been studied~\cite{AsllaniCarletti2018} assuming the connections between the oscillators to be given by a (weighted) all-to-all coupling, an Erd\H{o}s-R\'enyi network~\cite{ER,Bollobas} or a Newman-Watts network~\cite{NW}. In all the cases the resulting network is thus very dense and the variability in nodes connectivity is limited. In the present work we are interested in studying topologies where nodes exhibit an heterogeneous degree distribution and thus, besides the number of required controllers, a question that naturally arises is to determine their optimal ``location'' in the network. We thus consider and analyze the impact of different scenarios where controllers are selected uniformly at random among the available nodes, according to nodes degree or betweenness, and finally based on the recently introduced concept of functionability~\cite{RosellTarragoDiazGuilera}. The latter has been proposed as a single metric to rank nodes in the network according to their capacity to generate different behaviors once their dynamics is modified by some control term; interestingly enough the authors of~\cite{RosellTarragoDiazGuilera} have shown that more functional nodes have large degree and are peripheral ones.

By anticipating on the following, we will show that in case of Scale Free networks, controlling the hubs is a better solution than controlling any other node, to desynchronize the system. On the other hand, for networks with a strong core-periphery structure, we show that the best option is to control core nodes. It is interesting to observe that similar conclusions have been obtained in the case of synchronization \cite{WangChen2002,LuLiRong2009,Liu2018}. In those works, for instance, authors stated that if the number of pinned nodes is small, the best way to stabilize the synchronized state is a selection strategy based on the choice of the highest degree nodes, but this does not any longer hold true if the allowed number of controllers is large. We observe therefore a certain kind of universality of the efficiency of the selection strategies of pinned nodes, regardless of the actual goal of the control method.

The paper is organized as follows. In Section~\ref{sec:model} we present the Kuramoto model of coupled nonlinear oscillators, the control scheme proposed in~\cite{AsllaniCarletti2018} and a straightforward extension to allow to handle coupling structures closer to the ones presented in brain networks, i.e., with an heterogeneous degree distribution. Then in the following Section~\ref{sec:selmet} we introduce the selection schemes we have implemented to determine the number and the location of the controller nodes. Finally, after the presentation of the results in Section~\ref{sec:res}, we conclude and wrap up with Section~\ref{sec:conc}.

\section{The model}
\label{sec:model}

\begin{figure}
	\centering
	\includegraphics[width=0.325\linewidth]{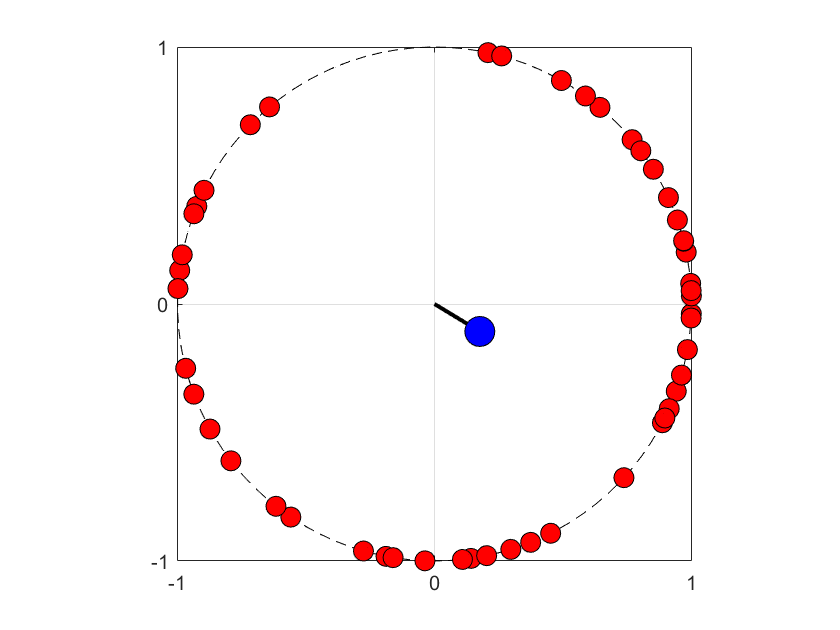}
	\includegraphics[width=0.325\linewidth]{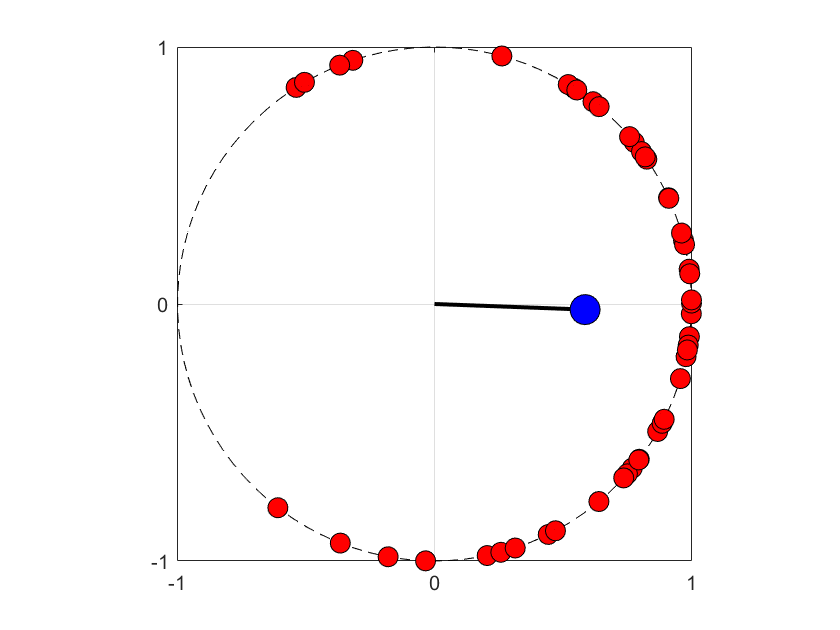}
	\includegraphics[width=0.325\linewidth]{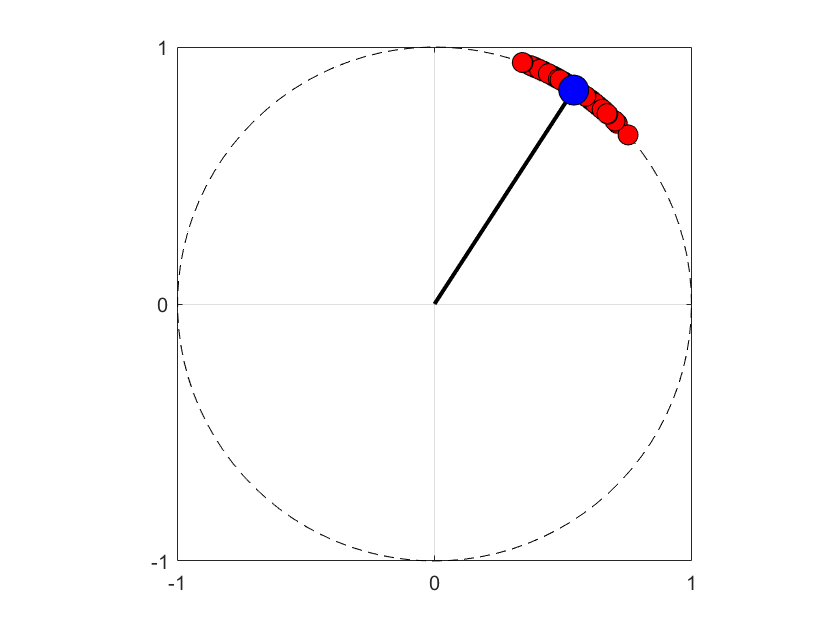}
	\caption{Representation of different synchronization states. Dots (red online) denote the Kuramoto oscillators at a fixed time, the (blue online) dot represents $R(t)e^{i\psi(t)}$ and the black segment its module. On the left panel the phases exhibit an asynchronous behavior and thus $R\approx0$; on the right panel the oscillators are almost perfectly synchronized and thus $R\approx1$. In the middle panel we report an intermediate case.}
	\label{Fig:order_parameter}
\end{figure}

To describe the behavior of basic units able to synchronize, representing, e.g., groups of neurons, we used the Kuramoto model (KM) of coupled oscillators that has been presented for the first time in~\cite{kuramoto1975}. Each one of the $N$ oscillators, is described by a phase $\phi_k$, with $k=1,\dots,N$, defined in $[0,2\pi]$ that evolves over time according to
\begin{equation} \label{eq:KM}
	\frac{d{\phi_k}}{dt}=\omega_k+\frac{K}{N}\sum_{j=1}^{N}A_{kj}\sin(\phi_j-\phi_k)\, ,
\end{equation}
where the first term on the right hand side, represents the oscillator natural frequency while the second one models the interaction existing among coupled oscillators. The latter term is defined by the adjacency matrix of the simple, undirected and unweighted network, $\mathbf{A}$,  namely $A_{ij}=1$ if the $i$-th and $j$-th oscillators do influence each other, and $A_{ij}=0$ otherwise. Finally, the parameter $K>0$ represents the coupling strength. Depending of the value of $K$, the distribution of the natural frequencies, $\Omega=(\omega_1,\dots,\omega_N)$, and the coupling topology, the KM can range from a synchronized state, where most of the units oscillate at unison in a coherent way, to an asynchronous one where on the opposite, almost each phase evolve independently from the other (See Fig.~\ref{Fig:order_parameter}).

The synchronization state of the system can be quantified by using the order parameter defined by
\begin{equation}
	\label{eq:Order_parameter}
	R(t)e^{i\psi(t)}:= \frac{1}{N}\sum_{k=1}^{N}e^{i\phi_k(t)}\, ,
\end{equation}
the latter represents thus the average of the phases, seen as arguments of complex numbers with unitary module. In case of synchronization, the values of the phases are close each other and thus they add each other in the sum defining $R$, the latter is thus large (see right panel of Fig.~\ref{Fig:order_parameter}). On the contrary in case of asynchronous behavior, the phases are homogeneously distributed in $[0,2\pi]$ and thus they sum to a small number, $R$ is then small (see left panel of Fig.~\ref{Fig:order_parameter}).

In some applications it is important to avoid a too strong synchronization state; scholars have thus designed control methods to reduce or suppress the synchronization in the KM. We hereby consider the framework proposed by~\cite{AsllaniCarletti2018} based on previous results~\cite{Witthaut} demonstrating that the KM dynamics can be embedded in the Hamiltonian system whose energy function is
\begin{equation}\label{eq::Hamiltionian}
H_0(\textbf{I},\Phi)=\sum_{i=1}^{N}\omega_iI_i-\frac{K}{N}\sum_{ij}A_{ij}\sqrt{I_jI_i}(I_j-I_i)\cos(\phi_j-\phi_i)\, , 
\end{equation}
where $\textbf{I}:=(I_1,\dots,I_N)$ is the action variable and $\Phi:=(\phi_1,\dots,\phi_N)$ the angle variable. More precisely, the torus
\begin{equation}\label{eq::torus}
	\mathcal{T}=\Big\{(\textbf{I},\Phi)|I_1=\dots=I_N=\frac{1}{2}\Big\}
\end{equation}
is invariant for the Hamiltonian flow and the dynamics on it is given by
\begin{eqnarray}
\begin{cases}\label{eq::HamOnTorus}
	\dfrac{dI_i}{dt}&=-\dfrac{\partial H_0}{\partial \phi_i}=0\\
	\dfrac{d\phi_i}{dt}&=\dfrac{\partial H_0}{\partial I_i}=\omega_i+\frac{K}{N}\sum_{j=1}^{N}A_{ij}\sin(\phi_j-\phi_i)
\end{cases}
\end{eqnarray}
that we can recognize to be the same ODE ruling the evolution of the KM~\eqref{eq:KM}. Authors of~\cite{Witthaut} have shown a relation between the synchronization of Kuramoto oscillators and the instability of orbits close to the torus $\mathcal{T}$. Starting from this observation and from ideas of Hamiltonian Control Theory~\cite{vittot,ciraolo}, it has been designed~\cite{gjata} a control term allowing to increase the stability of the torus $\mathcal{T}$ and thus to reduce the synchronization of the KM. More precisely one can modify Eq.~\eqref{eq:KM} into
\begin{equation} 
\label{eq:KMctrl}
	\frac{d{\phi_k}}{dt}=\omega_k+\frac{K}{N}\sum_{j=1}^{N}A_{kj}\sin(\phi_j-\phi_k)+S^{stim}_k(\phi_{i_1},\dots,\phi_{i_M})\, ,
\end{equation}
where $S^{stim}_k$ is a suitable control term~\cite{AsllaniCarletti2018} (see Appendix~\ref{app:control} for a more detailed description of the latter term). Let us emphasize that the control term depends on the $M\ll N$ controlled oscillators $\{i_1,\dots,i_M\}$. The rationale of the previous system~\eqref{eq:KMctrl} is that the signal is read from the controller oscillators, i.e., their phase is measured, then such information is elaborated and a new signal is computed and injected into the very same controller nodes and through them it reaches also their first neighbors. Interestingly such signal is always present but it is small, i.e., of order $K^2$, if the system is naturally in an asynchronous state while it increases as to compensate the KM behavior once the latter finds itself in a synchronous state. The control strategy is thus minimally invasive, the injected signal is often small and the number of controllers is limited $M\ll N$.

The aim of this work is to make one step further with respect to~\cite{AsllaniCarletti2018} along two research directions. First we analyze the behavior of the controlled KM on networks with heterogeneous degree distribution, in particular Scale Free networks and core-periphery ones, and second we propose a criterion to optimally select the location of the $M$ controllers that were randomly placed in~\cite{AsllaniCarletti2018}. Indeed, it seems quite natural to expect that acting on some nodes would have a greater impact on the system dynamics than other ones. A good selection criterion could therefore allow to decrease the number of required controllers and so minimize the control invasiveness by maintaining unchanged its efficiency.

\section{Selection methods}
\label{sec:selmet}

\begin{figure}[!t]
	\centering
	\includegraphics[width=0.31\linewidth]{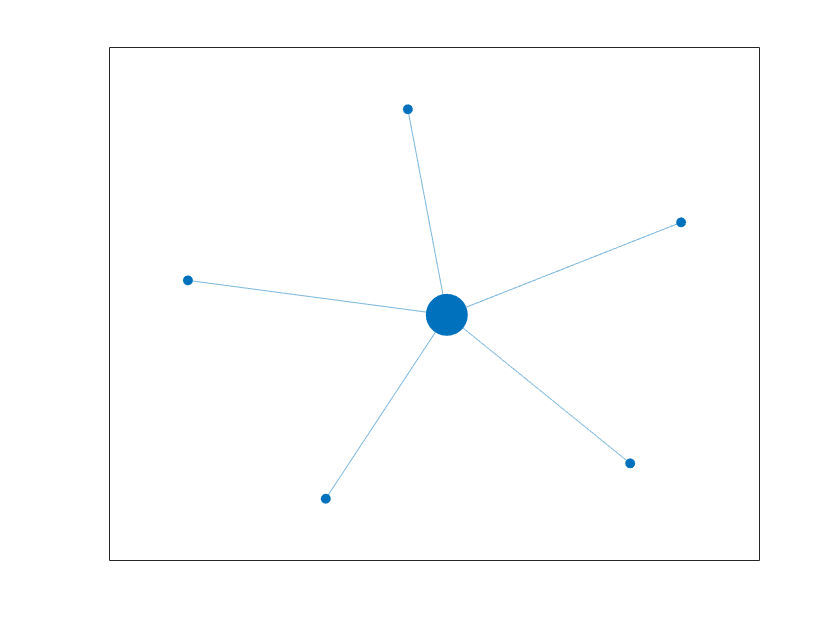}
	\includegraphics[width=0.31\linewidth]{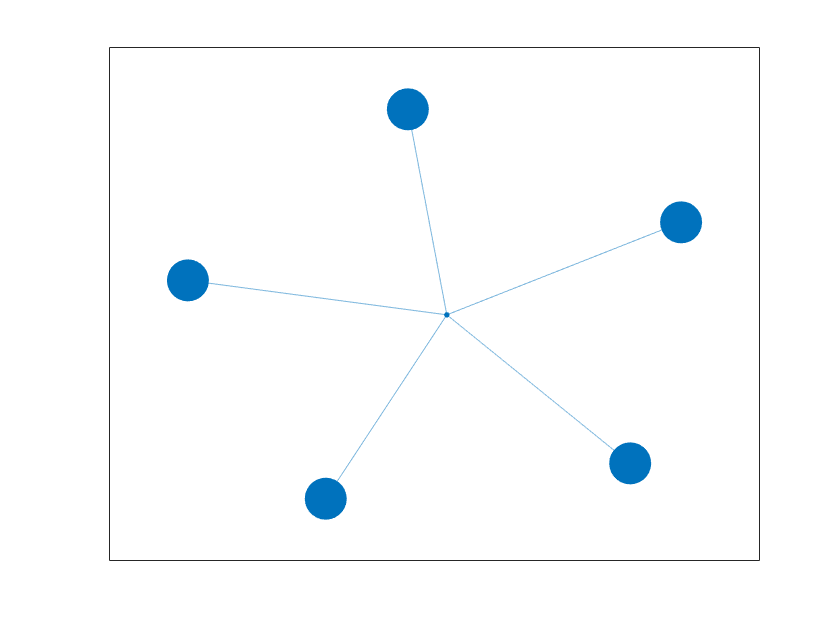}
	\includegraphics[width=0.31\linewidth]{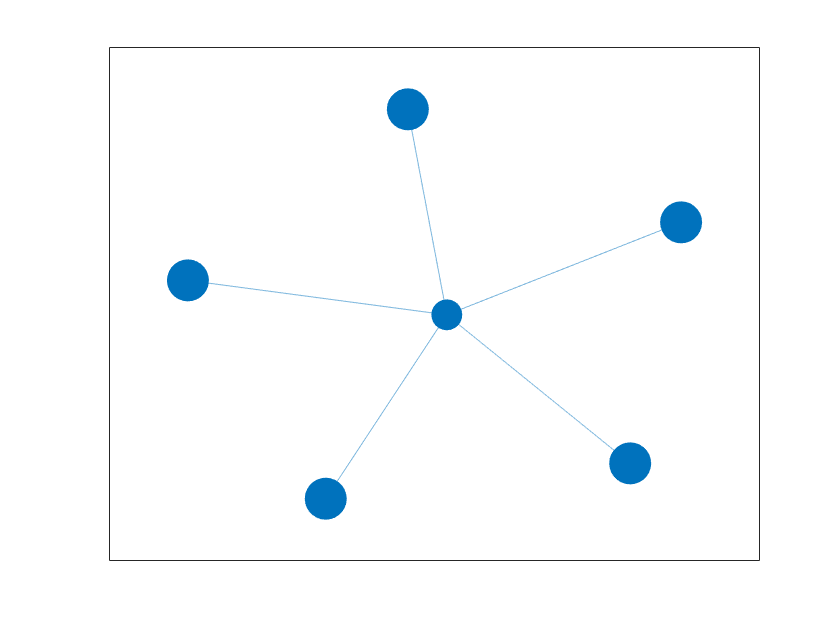}
	\includegraphics[width=0.31\linewidth]{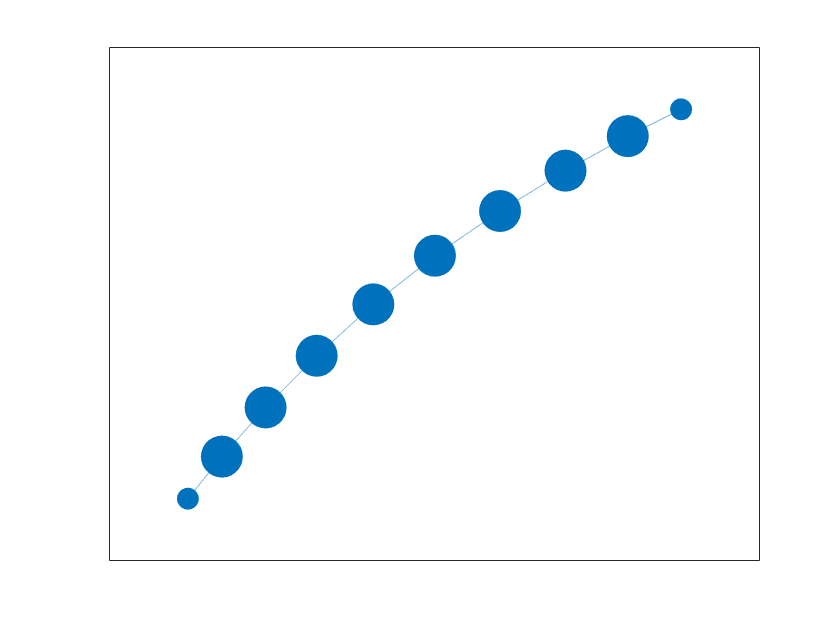}
	\includegraphics[width=0.31\linewidth]{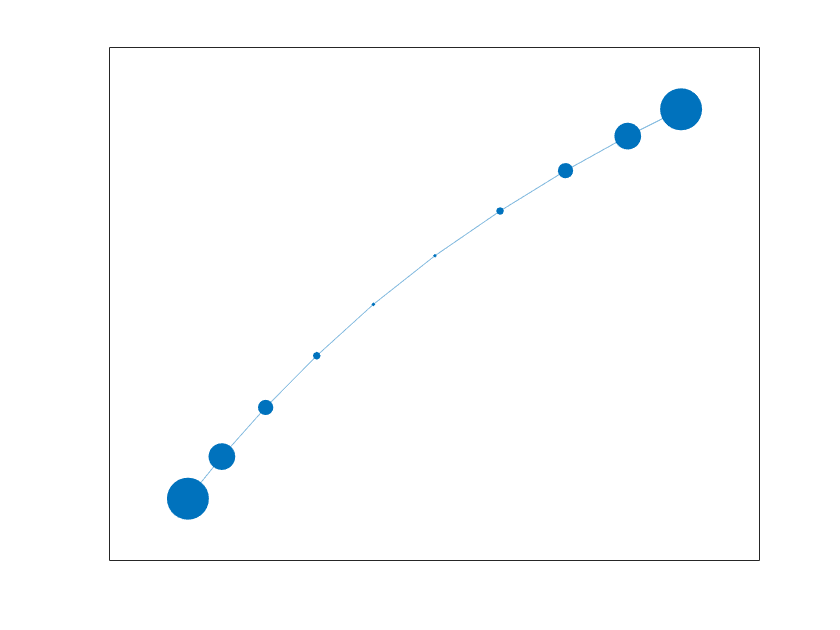}
	\includegraphics[width=0.31\linewidth]{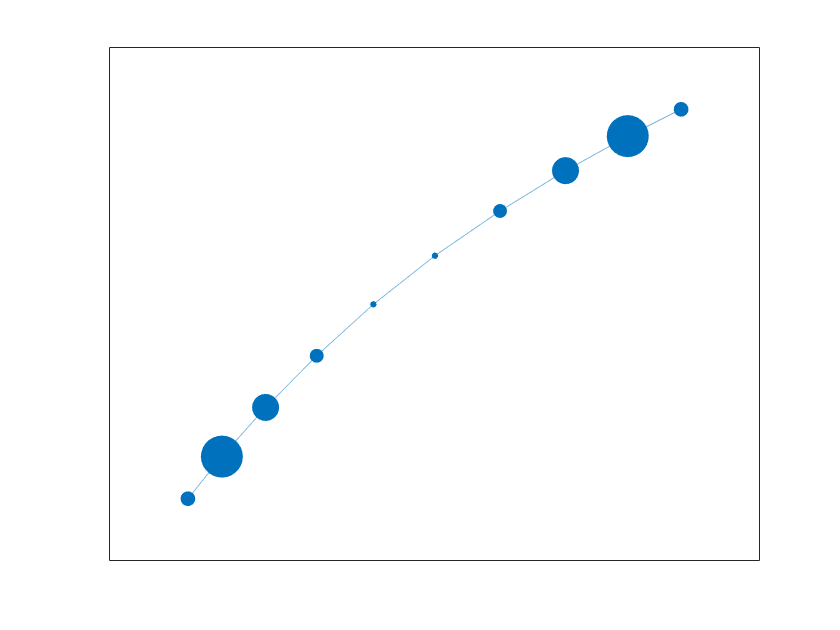}
	\caption{Compare node centrality, betweenness centrality and functionability. We show two simple networks where nodes size is proportional to their degree (left), the inverse of their betweenness (middle) and their functionability (right).}
	\label{Fig:functionability}
\end{figure}

As presented in the previous section, $M$ nodes should be selected to act as controllers; in this section we present the three used methods to perform such choice. In the first case we randomly select the $M$ nodes uniformly among the $N$ available ones; this method has been already used in~\cite{AsllaniCarletti2018} and proved satisfactory for networks without a structure, in the present case it would thus be used as a null-model and we refer to it as the {\em random selection}. The second choice is based on the node degree, namely the $M$ selected nodes are the $M$ nodes with the largest degree, hereby named {\em degree selection}. The last method is based on the functionability, a centrality measure developed in~\cite{RosellTarragoDiazGuilera}~\footnote{We also studied the selection method based on nodes betweenness, the results are reported in Appendix~\ref{app:betweenness}.}. Roughly speaking, it aims at measuring how influential is a node for the KM to reach a high synchronization state, thus $M$ nodes are chosen among those with the highest functionability score (see Appendix \ref{app:funct} for more details). This centrality measure is actually the result of a two antagonist competing features; indeed, it increases with the square of the node degree and, from its explicit formula, one can observe that it is anti-correlated with the betweenness centrality~\cite{RosellTarragoDiazGuilera}. Thus the functionability ranks high also peripheral nodes. In Fig.~\ref{Fig:functionability} we compare for two simple networks, a star and a line, the way nodes are ranked with respect to their degree (left panel), inverse of their betweenness (middle panel) and functionability (right panel). 

\section{Results}
\label{sec:res}

The aim of this section is to present the numerical results we realized both on synthetic networks (next section) and on few empirical ones (see Section~\ref{ssec:realnet}).

\subsection{Analysis of synthetic networks}
\label{ssec:syntnet}

\begin{figure}[!t]
	\centering
	\includegraphics[width=0.6\linewidth]{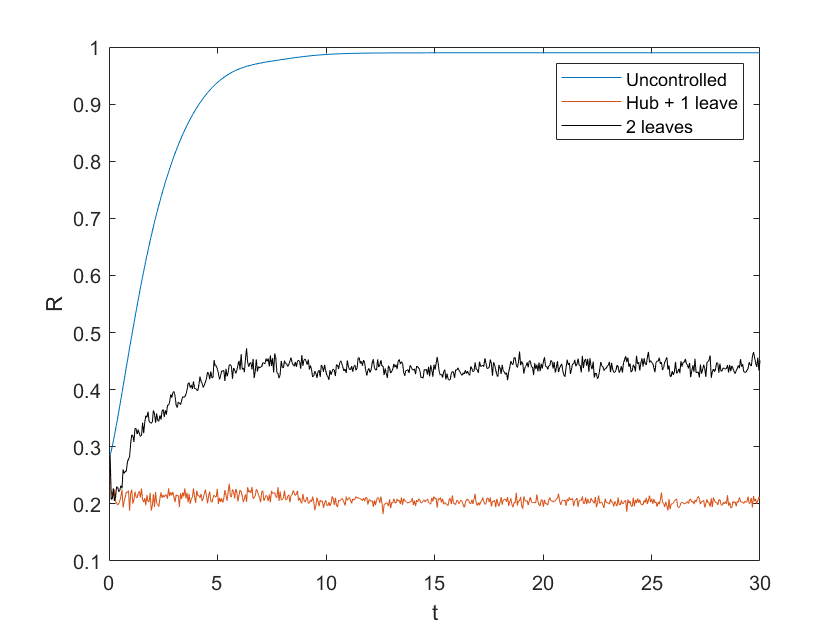}
	\caption{Evolution of $R(t)$ on a star network composed by $N=10$ oscillators modeled with the original KM, i.e., without controlled nodes (blue curve), with the controlled KM where the hub and one leaf have been pinned (red curve) and where two leaves act as controllers (black curve). Each curve has been obtained by averaging the results of $50$ independent simulations, where natural frequencies are distributed according to a normal distribution $\mathcal{N}(1,0.1)$ and initial conditions for angles are uniformly distributed in $[0,2\pi]$. We can observe that the original KM quickly synchronizes, the blue curve rapidly tends to a number close to $1$, while the controlled KM do not synchronize as testified by the low values of $R(t)$. Moreover the choice of controlling the hub seems to have the better desynchronization feature, indeed the red curve is well below the black one.}
	\label{Fig::star}
\end{figure}

Let us start to present our results applied to the simple and stylized star network shown in the top panels of Fig.~\ref{Fig:functionability}, where a central hub node is connected to $N-1$ leaves. Assume to have the possibility to use $M=2$ controllers, a natural question is thus to understand the behavior of the controlled system once one pinned node is the hub and the second one is chosen among the leaves or when both pinned nodes have been selected among the leaves. Which configuration returns better control features among the two? Let us observe that, although the natural frequencies are different from node to node, we considered all leaves as equivalent. Indeed, they are identical according to network topology, moreover averaging on several independent simulations obtained with different frequency, clearly smooth out possible differences. In Fig.~\ref{Fig::star} we report the results obtained for a star network made by $N=10$ nodes; we can observe that controlling the hub and one leaf (red curve) provides better results than controlling two leaves (black curve) because the order parameter $R$ remains quite low and far from the value $1$ reached by the uncontrolled system (blue line) for the whole considered time window.

Let us observe that similar conclusions have been obtained in the complementary framework of selecting the best nodes to achieve global synchronization by controlling a dynamical system with linear coupling (while KM exhibit a coupling through a sine function) and by using a linear control term. Authors have indeed demonstrated that one can control global dynamics of a scale free network by pinning the hubs~\cite{WangChen2002,Liu2018} and that in general the control efficiency depends on controllers degrees but also on the distance between them and the remaining nodes~\cite{LuLiRong2009,Liu2018}. 

On the other hand, as we have shown in Fig~\ref{Fig:functionability}, in a star network leaves are the nodes with highest functionability. So this simple example suggests that selection by functionability should work less efficiently than selection by node degree. This phenomenon will be confirmed by further observations.

Let us now consider the controlled KM to be defined on general scale free networks whose power law coefficient  $\gamma$ varies in  $[-2,-4]$. Such networks have been obtained by using the configuration method and the tools proposed in~\cite{LNR2017}. We generated $100$ networks~\footnote{Note that the configuration model do not guarantee that the resulting network is connected. In order to obtain connected networks of approximately $100$ nodes, we generated larger networks and extracted the maximal connected component whose size was between $85$ and $115$. So, not all the used scale free networks have the same size.} and we simulated the time evolution of the controlled KM defined on each network by using the three selection strategies: random, degree and functionability based. We then computed the asymptotic value of the order parameter $R_{\mathrm{as}}$, i.e., after a sufficiently long time, neglecting thus the transient phase, and averaging the order parameter over the remaining time series. Finally to compare the order parameters for networks of different sizes, we normalize the asymptotic order parameter by dividing it by the asymptotic value of order parameter obtained for the uncontrolled KM, $R$, obtaining in this way $\hat{R}=R_{\mathrm{as}}/R$.

\begin{figure}[!t]
	\includegraphics[width=0.33\linewidth]{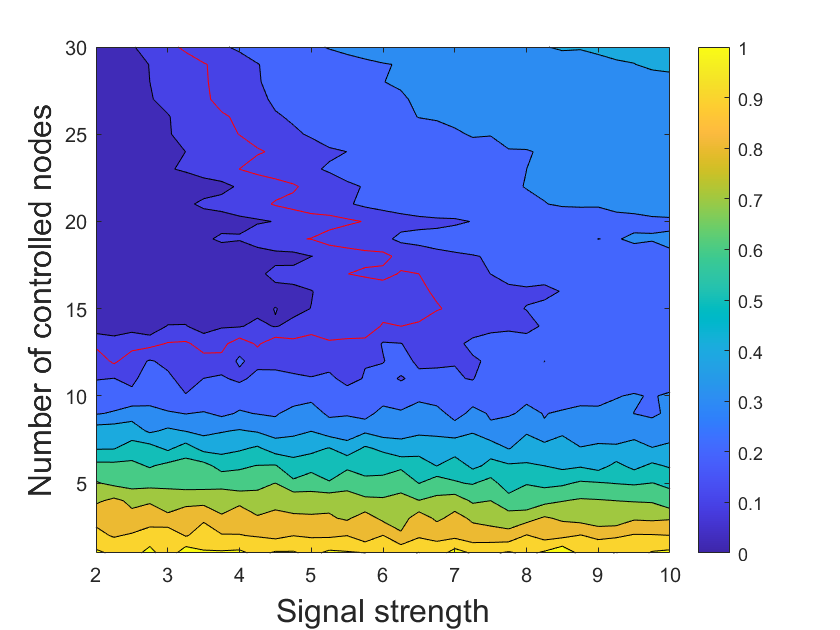}\hfill
	\includegraphics[width=0.33\linewidth]{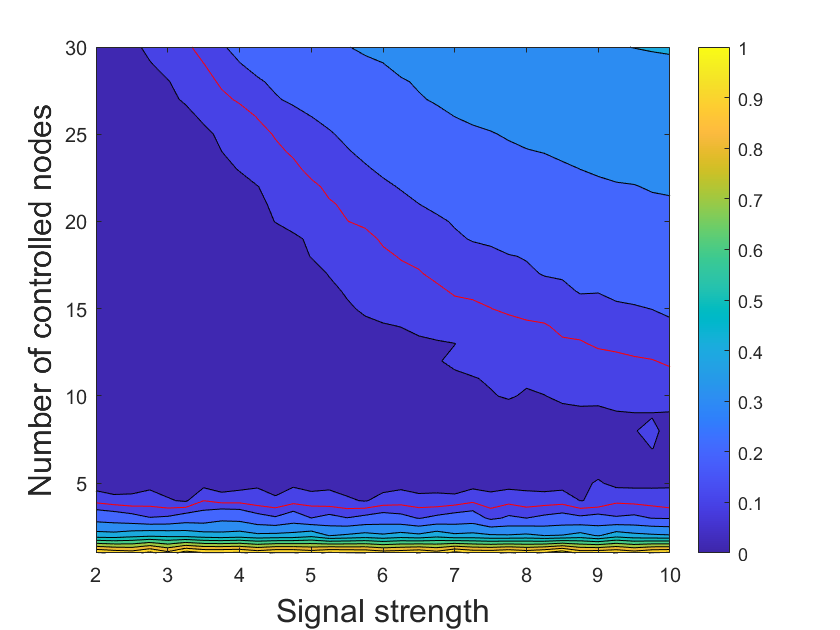}\hfill
	\includegraphics[width=0.33\linewidth]{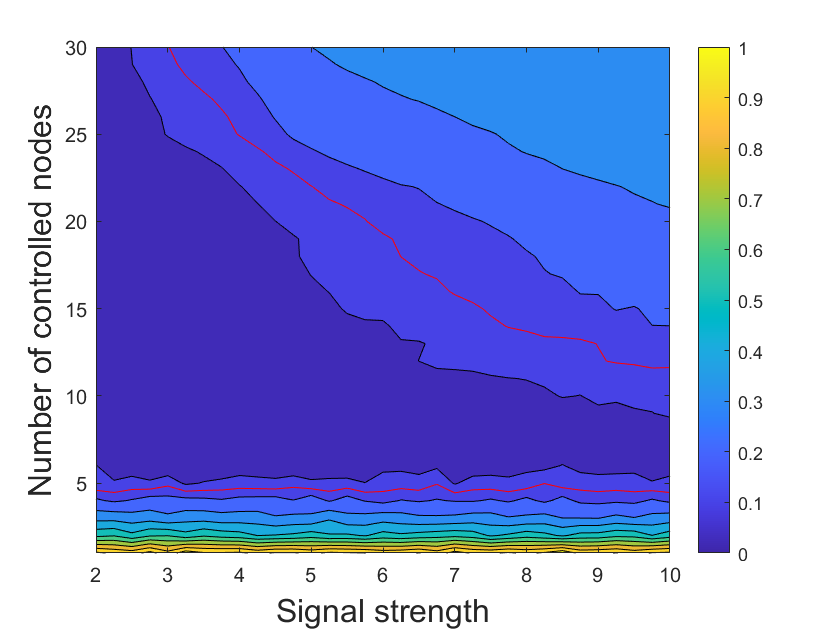}
	\caption{Normalized order parameter, $\hat{R}$, of the controlled KM defined on a scale free network, $\gamma=-3$, as a function of the number of controllers, $M$, and the control signal strength. By using a color code we show the value $\hat{R}$ averaged over $100$ scale free networks. Pinned nodes have been selected uniformly at random (left panel), on a degree-based selection (middle panel) and by using the functionability criterion (right panel). The red lines indicate the level contours at the threshold $\hat{R}=0.15$.}
	\label{Fig:sim_sf}
\end{figure}

\begin{figure}[t]
	\centering
	\includegraphics[width=0.6\textwidth]{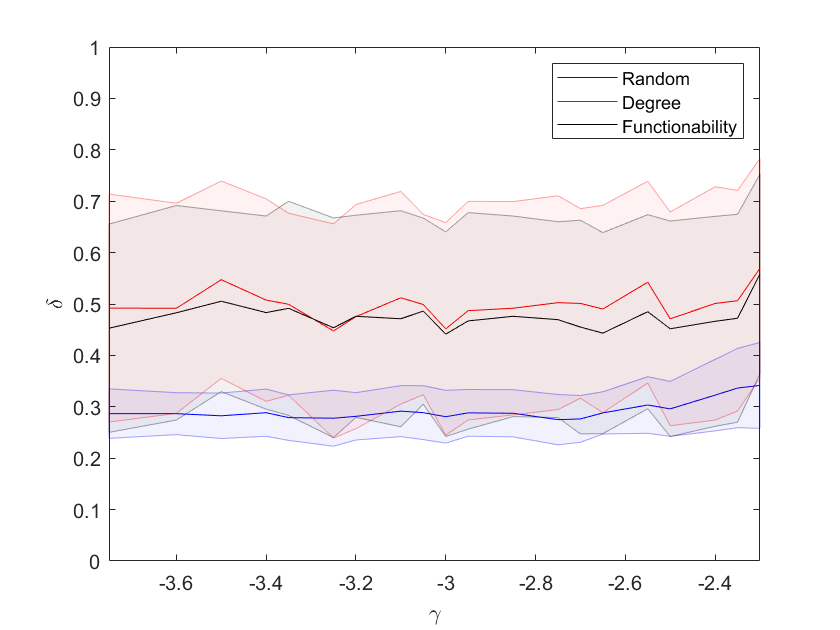}
	\caption{\label{Fig:sf_gamma} We present the size, $\delta$, of parameters region resulting into a well desynchronized state, $\hat{R}\leq 0.15$, as a function of $\gamma$ for scale free networks and the three used controllers selection methods. The darker lines correspond to the average of $\delta$ obtained by using $100$ independent networks realization, while the shaded areas about the average denote one standard deviation.}
\end{figure}

In Fig.~\ref{Fig:sim_sf} we report the average value of $\hat{R}$ computed from $100$ scale free networks with $\gamma=-3$ and in Table \ref{table:sim_sf} we report the fraction, $\delta$, of parameters values (number of controllers and control signal strength) associated to a well desynchronized state, for sake of definitiveness hereby defined by the condition $\hat{R}\leq 0.15$. From those results one can conclude that the degree selection allows to obtain better desynchronization results, indeed the parameters region corresponding to small value of $\hat{R}$ (dark blue) is larger in the case of degree selection than in the case of random or functionability ones and $\delta$ is the larger in this case. Similar results have been obtained for other values of $\gamma$ as we can see on Fig.~\ref{Fig:sf_gamma} where $\delta$ is plotted as function of $\gamma$ for the three selection strategies. Indeed, one can appreciate that for all value of $\gamma$ taken into account the degree selection always shows (slightly) better results than functionability, and both do far better than the random one. Moreover, we can observe than the efficiency of the control seems to be independent of the exact value of $\gamma\in[-2,-4]$. The curves presented in Fig.~\ref{Fig:sf_gamma} do not show any particular trend as functions of $\gamma$. So it appears that, as long as the degree distribution is assumed to be heterogeneous (as in the case of a power-law degree distribution with a negative exponent), the exact distribution parameters (here $\gamma)$ does not affect the efficiency of the degree selection strategy. However, we can notice that even if degree- and functionability-based strategies generally produce better results than the random selection, the variances of results for the formers are far more large. 

\begin{table}[t]
	\centering
	\begin{tabular}{|c|c|c|c|}
		\hline
		&\textbf{Random selection}&\textbf{Degree selection}&\textbf{Functionability selection}\\
		\hline
		$\delta$&0.2091&0.5626&0.5141\\
		\hline
	\end{tabular}
	\caption{\label{table:sim_sf} We report the fraction $\delta$ of parameters, number of controllers and control signal strength, for which $\hat{R}\leq 0.15$ for the controlled KM defined on a scale-free networks with $\gamma=-3$. These values quantify the level sets of $\hat{R}$ shown with a color code in Fig.~\ref{Fig:sim_sf}.}
\end{table}

The numerical results support the claim that degree-based and functionability-based selection strategies both provide a stronger reduction of the synchronization with respect to a random selection of the pinned nodes. Indeed, low values of the averaged $\hat{R}$ (darker blue area in Fig.~\ref{Fig:sim_sf}) can be reached already with $5$ controllers chosen by their degree or functionability (middle and right panels), while this is not the case with a random selection (left panel). The similarity of the results obtained with those two methods can be justified by the fact that quite often the same nodes are selected by the two strategies. Indeed, for the scale free networks if we compare the values of the degree and the functionability of each node, the mean correlation coefficient between those two measures is $\approx0.77$. Thus, hubs are also often the nodes with highest functionability score and the two rankings used for the nodes selection are very similar. In this case, it thus seems that the dependence of the functionability score on the inverse of betweenness is somehow reduced with respect to the node degree, possibly because of the huge heterogeneity of the degree distribution in scale free networks. Based on this observation and on the fact that selection by degree is computationally cheaper, we would prefer the latter strategy selection. Let us recall that the functionability also promotes nodes with low betweenness, however (see Appendix~\ref{app:betweenness}) our results support the claim that this factor does not have a strong influence on the desynchronization.

Let us conclude by observing a somehow counter-intuitive phenomenon in the case of degree and functionability selection. By using a too large number of controllers and a too large signal strength does not help to reduce the synchronization; indeed values of $\hat{R}$ associated to parameters lying in the the top-right regions (see Fig.~\ref{Fig:sim_sf}) are larger than the values obtained by using a smaller number of controllers or a smaller control signal strength. We can interestingly draw a parallel between those latter observations and other results available in literature about pinning control methods developed to achieve synchronization. For instance in~\cite{Liu2018}, authors concluded that the best controllers subset is composed preferentially by the hubs when the number of controllers is small but this is not the case when the number of pinned nodes is large. It seems that the conclusions are the same in the presented framework as one can observe from Fig.~\ref{Fig:sim_sf} where values for $\hat{R}$ associated to $5<M<10$ are smaller than those computed for $25<M<30$, and this despite the fact that the goal of the control is not the same: we want to achieve desynchronization and not enhance synchronization.

To conclude this section let us present the numerical results obtained by studying the KM defined on core-periphery networks. The latter are characterized by a relatively small subset of nodes densely connected among them, forming thus a core, and the presence of peripheral nodes, weakly connected among them while strongly connected to the core nodes~\cite{Borgatti}. The network exhibits thus an heterogeneous degree distribution as core nodes have many connections by contrast with periphery ones, but the number of hubs is in general larger than in scale free networks. Moreover, such network exhibits the rich club property~\cite{ZhouMondragon,Colizza}; because most of the information passes through the core such networks are thus a generalization of the simple star network previously discussed. Let us observe that a similar structure with central node can also be found in the brain~\cite{vdheuvel}.

To generate core-periphery networks we start by building a complete network made of $n<N$ nodes, that will eventually be the core. Then, periphery nodes are successively added one by one and each of them sequentially linked to one core node with probability $p$, selected uniformly at random, and with probability $1-p$ to one non-core node, also selected uniformly at random. After having added sufficiently many nodes we can consider the periphery structure to have been built. The process gives rise to a network with a star shape but with a large core of $n$ nodes instead of a single hub, and with several branches of different lengths, instead of the leaves, that are formed by the periphery nodes. 

We are then interested in studying the controlled KM defined on top of such core-periphery network and in particular to understand the impact of the location of the controllers chosen among the core nodes or the periphery ones. For this reason we fix the number of available controllers, $M$, and we consider $k$ pinned nodes, $k=0,\dots,M$, to belong to the core, and thus $M-k$ peripheral nodes.

\begin{figure}[ht]
	\centering
	\includegraphics[width=0.41\linewidth]{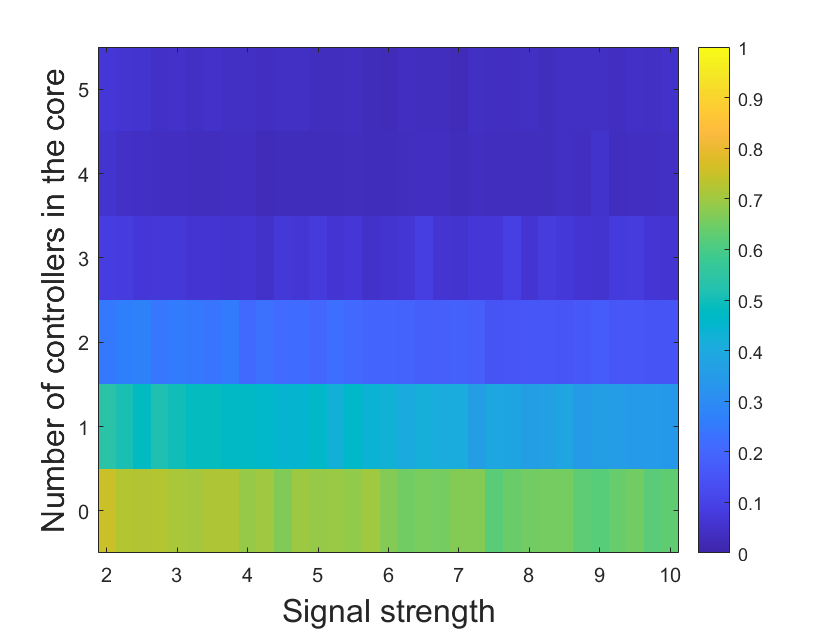}\quad
	\includegraphics[width=0.41\linewidth]{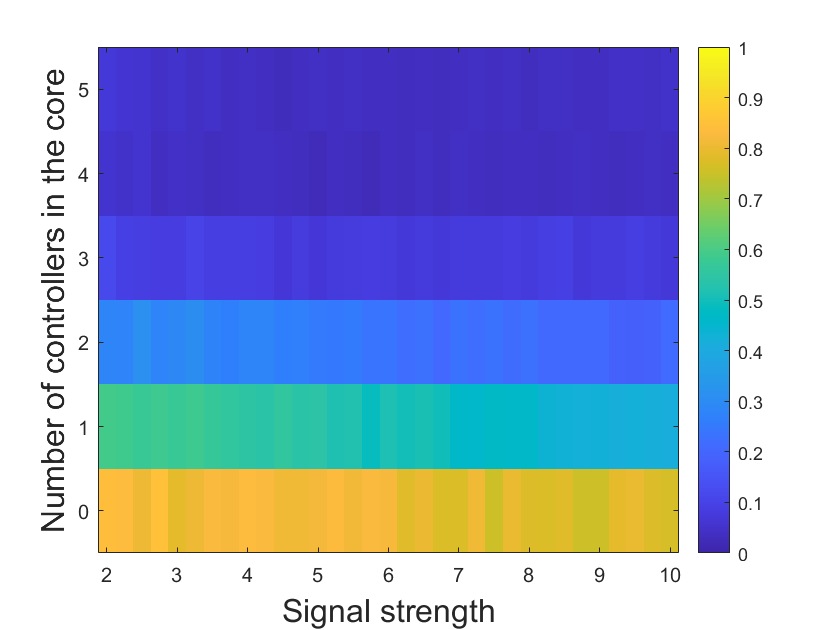}
	\caption{
	We show the normalized order parameter computed for controlled KM defined on core-periphery network, with $n=10$ nodes in the core and with $M=5$ controllers, as a function of the number of controllers in the core, $k$, and the control signal strength. By using a color code we show the value $\hat{R}$ averaged over $100$ core-periphery networks with $p=0.7$. Left panel corresponds to the node degree-based selection, while the right one to the functionability-based strategy.}
	\label{Fig::cp10} 
\end{figure}

For sake of definitiveness let us start to deal with a core of size $n=10$, $p=0.7$ and consider $M=5$ controllers. Then for $k=0,\dots,5$ we select $k$ nodes in the core according to their degree or to their functionability score (i.e., the $k$ core nodes with the highest score for the considered centrality measure) and we do the same with $5-k$ nodes in the periphery. We then numerically solve the controlled KM. Once the orbits have been obtained we compute as before the normalized order parameter $\hat{R}$ and its average over $100$ independent replicas (corresponding to the number of different core-periphery networks generated as explained above). Results reported in Fig.~\ref{Fig::cp10} show that the more there are controllers in the core, the smaller is the averaged $\hat{R}$, hence the better the desynchronization process. For instance controlling $3$ core nodes selected by degree or functionability guarantees to substantially reduce the synchronization with respect to the case where the core hosts $2$, $1$ or $0$ controllers. A similar trend is observed also for uniformly random selection (see Fig.~\ref{Fig::app::cp_p} in Appendix~\ref{sec:figureapp}) but one need to impose more controllers in the core (at least $4$) to fully desynchronize the system.

As with the previous case of scale free networks, degree and functionability selections provide similar results and, still, the former one seems to be a slightly better criterion. Indeed, the averaged $\hat{R}$ values obtained are lower with the degree-based method (left panel) than with the functionability-based selection (right panel). The results do agree with the previous ones obtained by using scale free network because they point out the relevance to pin nodes with high degree to reduce synchronization. Let us observe that similar results can be obtained by using a larger core and more controllers (see Appendix~\ref{sec:figureapp} for the case $n=20$ and $M=7$).

\begin{figure}
	\centering
	\includegraphics[width=0.75\linewidth]{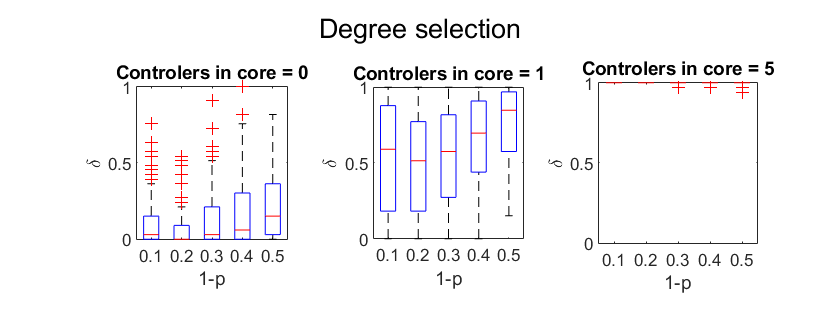}
	\includegraphics[width=0.75\linewidth]{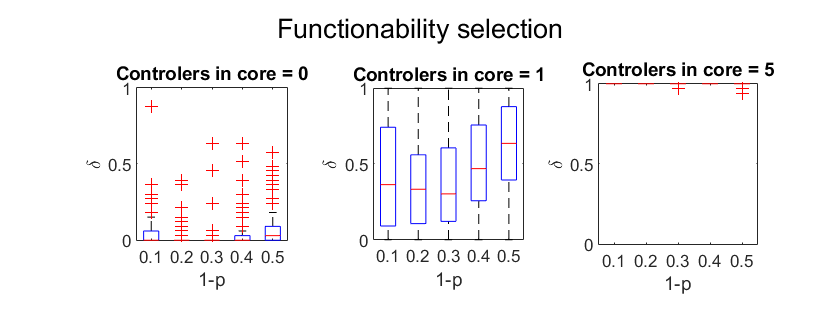}
	\caption{We show the distribution of the size, $\delta$, of parameters region resulting into a well desynchronized state, i.e., associated to $\hat{R}\leq 0.15$, for different core-periphery networks made by $N=100$ nodes, generated for several values of $p$. Top panels: pinned nodes are selected according to their degree; bottom panels: selection is performed according to functionability. Left columns correspond to $0$ controlled nodes in the core, middle columns correspond to $1$ pinned node in the core and $4$ in the periphery, while the right columns show the cases where $5$ pinned nodes are in the core.}
	\label{Fig::cp_p} 
\end{figure}

To finish this part about core-periphery networks, let us study the efficiency of the control and the selection strategies in function of the parameter $p$ that determine the size of the periphery. Because we want that with large probability a new node is added to the core we consider $p\in[0.5,1]$, this allows thus to avoid cases where the periphery nodes become more ``central'' than core nodes, i.e., with larger degree. For any given $p$ we computed again the normalized order parameter $\hat{R}$ and the size of the parameters region, $\delta$, for which $\hat{R}\leq0.15$, i.e., the fraction of dark blue area in synchronization maps. The results are shown in Fig.~\ref{Fig::cp_p} and we can again draw the general conclusion that putting a large number of controllers in the core is the best strategy, independently on the way they will be exactly chosen, i.e., degree or functionability, and on $p$. Moreover, we observe that once some controllers are in the periphery, the performances are better when $p$ is close to $0.5$. This can be easily explained by observing that with $p\sim 0.5$, periphery nodes are densely connected each other and less with the core nodes and, therefore, it is advantageous to choose some of them as controllers.

All these observations suggest that, even if controlling high degree nodes is important, to spread controllers in the network can have as well a relevant impact. Indeed, pinned core nodes are close each other (once we measure the distances in terms of links separating two nodes) and thus they generally lower the average distance between the controllers and the unpinned nodes. Moreover, in particular configuration in which $p\sim 0.5$, choosing some periphery nodes as pinned ones can be appropriate because a controller in the ``middle'' of the tail would reduce the distances between the controlled and the uncontrolled set of nodes.

Before moving to the conclusions, notice that we have also analyzed the above selection strategies on small-world networks obtained by using the Watts-Strogatz algorithm~\cite{WS1998}. In this case, the degree distribution is far more homogeneous and both selection strategies are not really satisfying, even if the degree-based one provides quite good results. This behavior can be explained by recurring again on the relevance of an optimal spreading of the controllers among the network nodes, i.e., the average distances between controllers and other nodes has to be small (see Appendix~\ref{app:ws} for more details). 

\subsection{Real networks}
\label{ssec:realnet}
Let us now consider the desynchronization capabilities of the control strategies applied on real empirical networks. We decided to use three networks taken from~\cite{networkswebsite}, representing brain connectivity of, respectively, cat, macaque and mouse. They were selected because of their very first nature, i.e., brain networks, and also because they have small enough sizes to be able to obtain good statistics in a reasonable time, being their sizes respectively: $N=65$, $N=91$ and $N=213$. Let us observe that each network has been transformed into a symmetric one to fit with our working assumption.

\begin{figure}[!h]
	\includegraphics[width=0.33\linewidth]{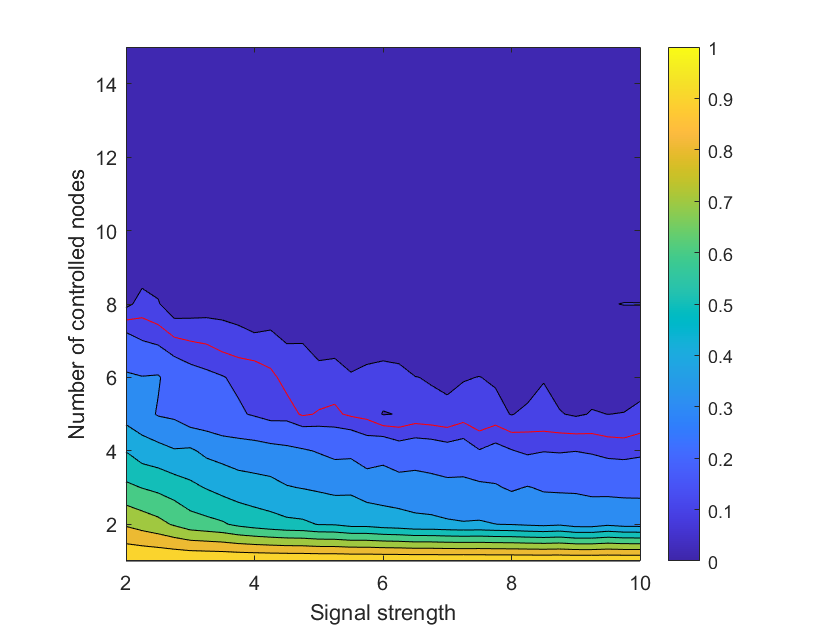}\hfill
	\includegraphics[width=0.33\linewidth]{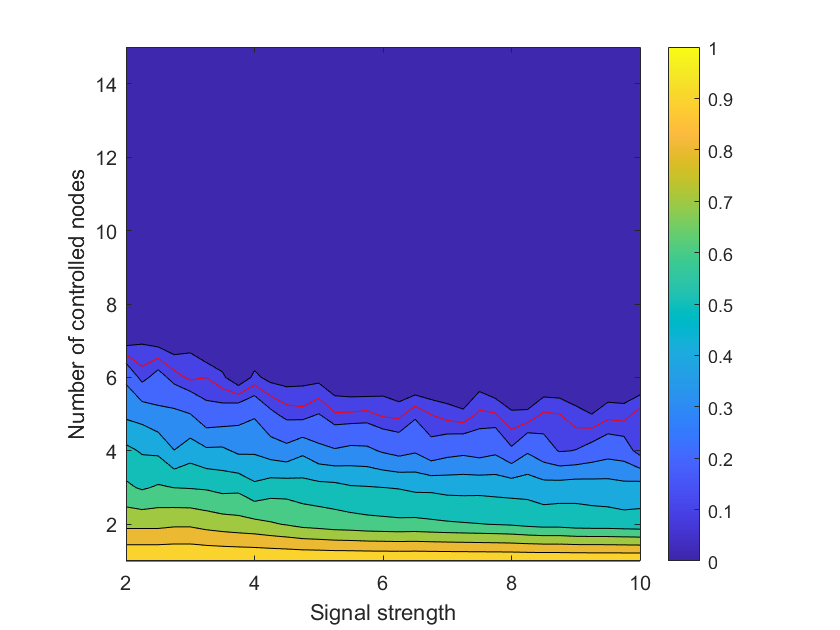}\hfill
	\includegraphics[width=0.33\linewidth]{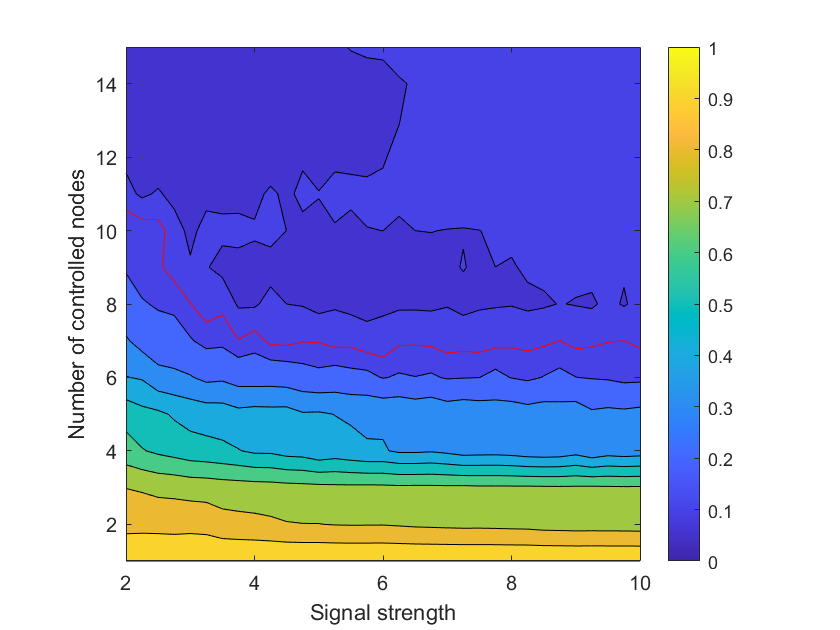}\\
	\includegraphics[width=0.33\linewidth]{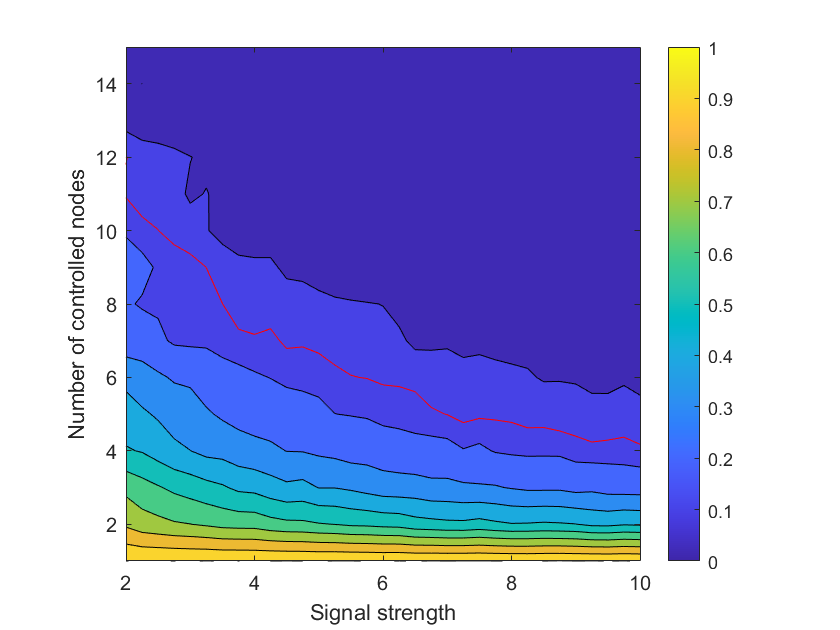}\hfill
	\includegraphics[width=0.33\linewidth]{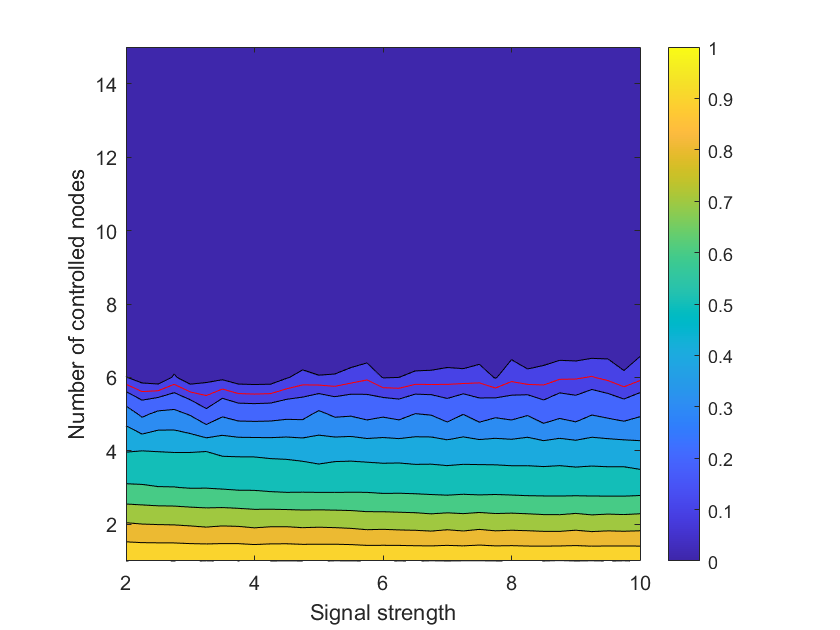}\hfill
	\includegraphics[width=0.33\linewidth]{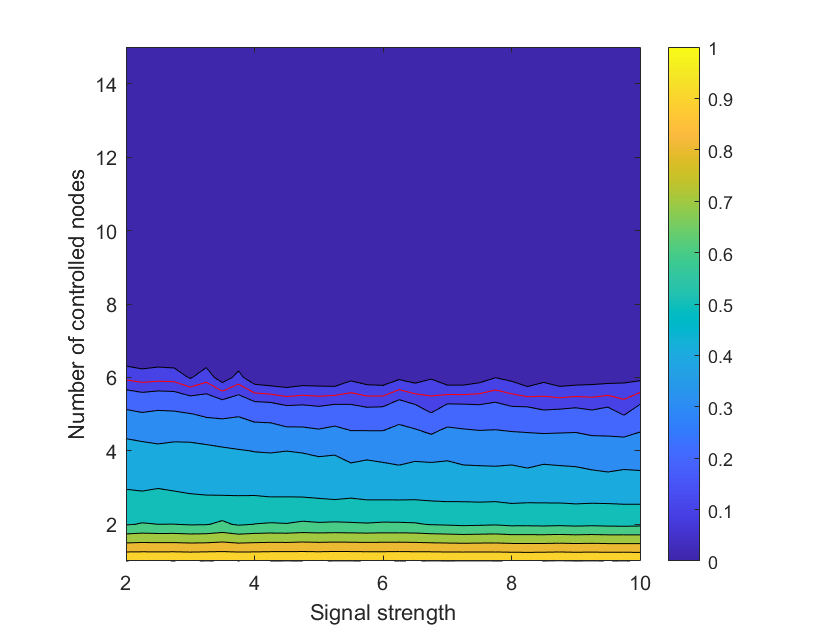}\\
	\includegraphics[width=0.33\linewidth]{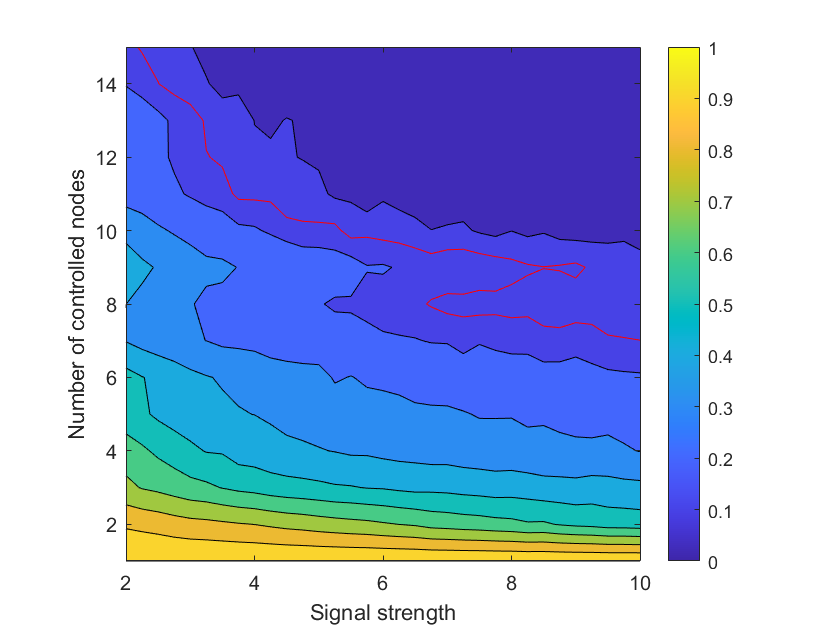}\hfill
	\includegraphics[width=0.33\linewidth]{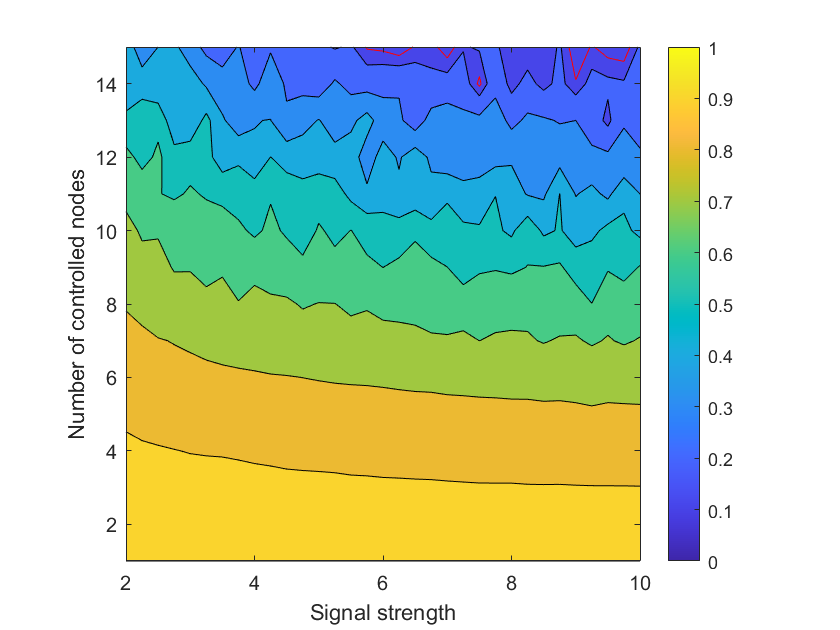}\hfill
	\includegraphics[width=0.33\linewidth]{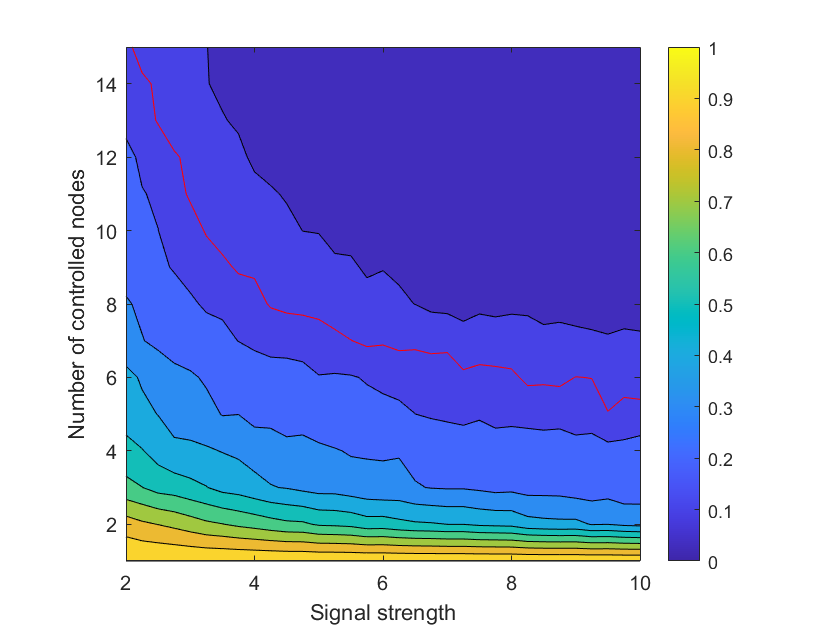}
	\caption{We show the normalized order parameter $\hat{R}$ obtained from numerical simulations of the controlled KM applied on the three real networks extracted from~\cite{networkswebsite}. For each couple, control signal strength and number of controllers, $50$ simulations have been repeated with a fixed network and randomly assign natural frequencies $\omega$ following a normal distribution $\mathcal{N}(1,\sigma^2)$, $\sigma=0.1$. The first row corresponds to cat's connectome, the second one to the macaque and the third one to the mouse. The first column refers to the random selection, the second to the degree one and the third to functionability. The red lines indicate the contour level corresponding to $\hat{R}=0.15$.}
	\label{fig::real_networks}
\end{figure}

The $\hat{R}$ values obtained for each network and each selection strategy by using the same ideas as before, are presented in Fig.~\ref{fig::real_networks} as function of $M$ and the control signal strength; in Table~\ref{table:sim_real} we report the fraction $\delta$ of parameters values for which $\hat{R}\leq 0.15$. We  can generally observe the same kind of results than the ones described above, namely the degree selection strategy provides better results, except for the mouse connectome where the degree-based selection seems to be less efficient than both other strategies. 

\begin{figure}
	\includegraphics[width=0.3\linewidth]{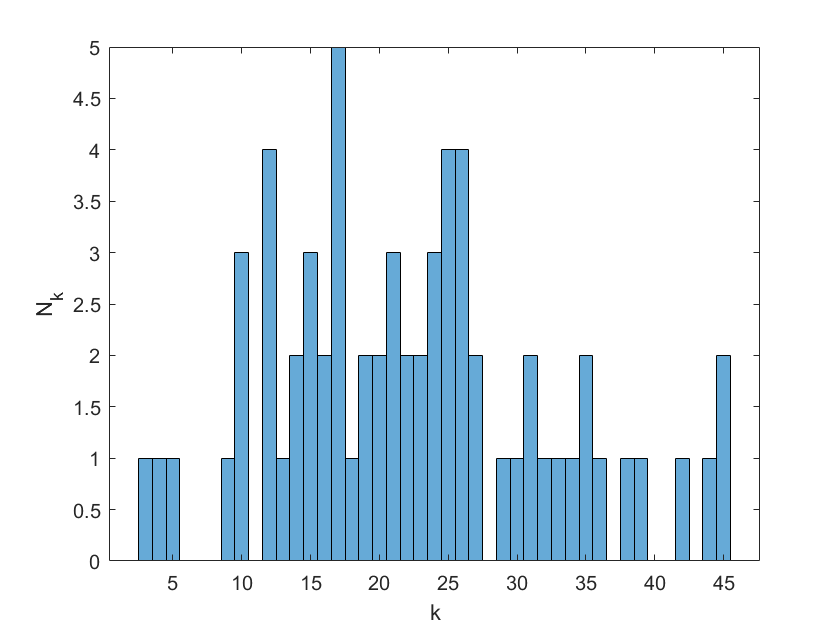}
	\includegraphics[width=0.3\linewidth]{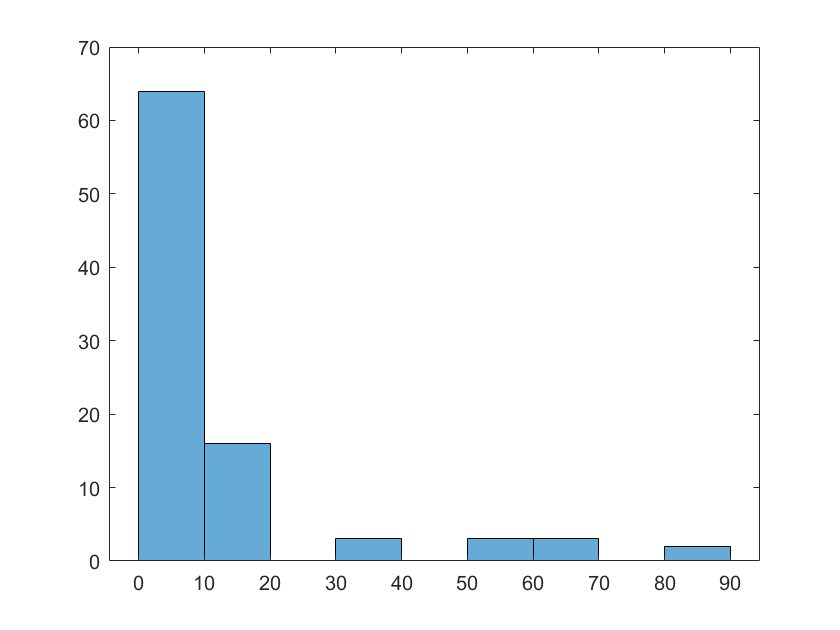}
	\includegraphics[width=0.3\linewidth]{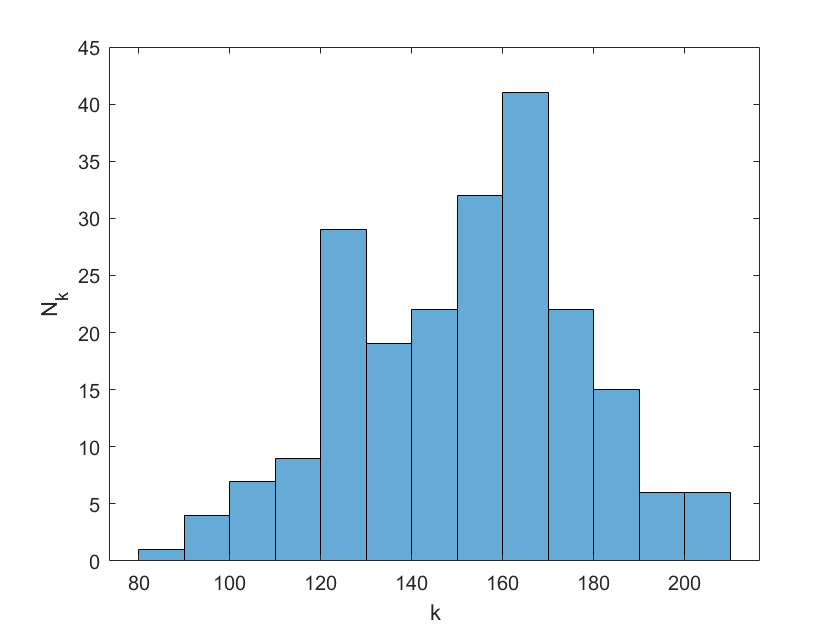}
	\caption{\label{fig::real_net_degree_dist} Degree distributions of three real connectomes taken from~\cite{networkswebsite} (number $N_k$ of nodes with degree $k$). Left panel: cat connectome; middle panel: macaque connectome; right panel: mouse connectome.}
\end{figure}

An explanation of the latter phenomenon can be found by analyzing the degree distribution of the three graphs (see Fig.~\ref{fig::real_net_degree_dist}). In the mouse connectome case, the network is very densely connected and the number of high degree nodes is very large. A vast majority of nodes are connected to more than $N/2$ neighbors. In this context, we can draw a parallel with previous observations about scale free networks: pinning a too large number of hubs is not a good strategy because the amount of control signal injected in the system is too strong and applied to too much close nodes, inducing thus synchronization of the latter. Indeed, the control term can become dominant and create a new kind of system trajectory around which the system can synchronize. For this reason selecting nodes by functionability, or even randomly, can better induce more subtle perturbations and by consequence brings the system to a desynchronized state. 
\begin{table}[ht]
	\centering
	\begin{tabular}{|c|c|c|c|c|}
		\hline
		&&\textbf{Random selection}&\textbf{Degree selection}&\textbf{Functionability selection}\\
		\hline
		{\textbf{Cat}}
   		&\textbf{Original}&0.6788    &0.6808    &0.5556\\
		\hline
		{\textbf{Macaque}}
		&\textbf{Original}&0.6101    &0.6646    &0.6667\\
		\hline
		{\textbf{Mouse}}
		&\textbf{Original}&0.3636    &0.0162    &0.5152\\
		\hline
	\end{tabular}
	\caption{
	We report the fraction $\delta$ of parameters, number of controllers and control signal strength, for which $\hat{R}\leq 0.15$ computed for the controlled KM defined on the three real connectome networks, such values correspond to the Fig.~\ref{fig::real_networks}.}
	\label{table:sim_real}
\end{table}

In Appendix~\ref{sec:reshuffled} we report results obtained by using degree-preserving reshuffled versions of those three networks, interestingly those results are very close to the ones for the original network, supporting thus the claim that the degree distribution is a key factor to determine the desynchronization feature of the controlled KM.
In addition we can note that, in the case of cat and macaque networks, where the number of hubs is relatively low (see Fig. \ref{fig::real_net_degree_dist}), the results are aligned with the previous observations. 

\section{Conclusions}
\label{sec:conc}

In this work we have studied the possibility to desynchronize a network of coupled Kuramoto oscillators by using a non-invasive control strategy developed in~\cite{gjata,AsllaniCarletti2018}. By using synthetic networks and empirical brain connectomes we have studied the impact of the location of the controllers and their number in the desynchonization efficiency, i.e., the reduction of the Kuramoto order parameter with respect to the uncontrolled case. The results of this work improve thus over the ones presented in~\cite{AsllaniCarletti2018} along two direction: first, because we consider heterogeneous networks, synthetic and real ones, second we consider several selection strategies for the controlled nodes. As a general conclusion, we can observe that our results align with previous works aiming at controlling linearly coupled system with the goal of increasing the synchronization~\cite{WangChen2002,LuLiRong2009,Liu2018}, hence in a complementary framework than ours. That sheds light on the existence of a  possible universal rule to achieve optimal node selection in the field on network control theory. Indeed, we found that when there is a small subset of high degree nodes as in scale free or core-periphery networks, the optimal way to reduce synchronization is to inject the control signal in the hubs. The degree-based selection results thus to be optimal in scale free networks and, in general, concentrating the control on the core in core-periphery networks is more efficient than periphery pinning. 

On the other hand, another important feature is minimizing the mean distances between the controllers and the other nodes. Indeed, this claim is supported by our study of the core-periphery networks when $p$ is small and also with small-world networks, where the control efficiency is anti-correlated with the mean distances between the controllers and the other nodes and correlated with the minimal (positive) eigenvalue of the reduced Laplacian matrix $\lambda_{min}(\tilde{L})$ (see Appendix~\ref{app:ws}). As an illustration, functionability-based selection is not efficient on small-world networks and gives rise to low $\lambda_{min}(\tilde{L})$ and high mean distances values.

According to these observations the degree-based selection seems to be in general the best way to select controllers, even if the degree distribution is quite homogeneous. Indeed, it always provides at least slightly better results than functionability selection, and random one, and it is less computationally demanding than the former. Moreover, especially if the degree distribution is heterogeneous, it is important to control the dynamics of the nodes with many connections to better spread the control effect, because direct action of the pinned nodes onto their first neighbors, mimicking the action of the electromagnetic field effect.

Let us observe that some questions still remain open. For example, in core-periphery networks, by using random selection, one can observe that in certain cases controlling some periphery nodes is also favorable; this suggests that other node features than degree have to be considered in the selection to improve control efficiency. It is also worth mentioning the loss of efficiency in desynchronization in the case of scale free networks endowed with many controllers and large signal strengths, our intuition is that in this way clusters of strongly synchronized nodes are formed; let us observe that this this fact strengthens the parallel with previous results devoted to determine control methods to enhance synchronization~\cite{Liu2018} whose conclusion is that pinning the hubs is not the best strategy in the case of a large number of controllers. Finally it would be interesting to consider ``spatial'' networks~\cite{Barthelemy} where, hence, distances can be computed by using, e.g., the Euclidean distance, in this way the impact of pinned nodes on their neighbors would be driven by geometrical positions and not number of links. We believe that these issues are relevant but they go beyond the scope of this work and we left it thus to a further study.

\section*{Acknowledgements}

Part of the results were obtained using the computational resources provided by the "Consortium des Equipements de Calcul Intensif" (CECI), funded by the Fonds de la Recherche Scientifique de Belgique (FRS-FNRS) under Grant No. 2.5020.11 and by the Walloon Region.

\appendix

\section{Derivation of the control term}
\label{app:control}

It has been shown~\cite{AsllaniCarletti2018,gjata} that the Kuramoto Model~\eqref{eq:KM} can be embedded in the Hamiltonian system
\begin{equation}\label{eq::derHam}
\begin{cases}
\dot{I_i}=-\dfrac{\partial H_0}{\partial \phi_i}=-2\frac{K}{N}\sum_{j=1}^{N}A_{ij}\sqrt{I_iI_j}(I_i-I_j)\cos(\phi_j-\phi_i)\\
\dot{\phi_i}=\dfrac{\partial H_0}{\partial I_i}=\omega_i+\frac{K}{N}\sum_{j=1}^{N}A_{ij}\left[2\sqrt{I_iI_j}\sin(\phi_j-\phi_i)-\sqrt{I_j/I_i}(I_j-I_i)\sin(\phi_j-\phi_i)\right]
\end{cases}\, ,
\end{equation}
obtained from the Hamiltonian function~\eqref{eq::Hamiltionian}. Moreover, the torus $\mathcal{T}$, given by~\eqref{eq::torus}, is invariant and the Hamiltonian flow restricted on it, corresponds to the Kuramoto system. Finally, in~\cite{Witthaut} authors shown the existence of a links among the  stability property of $\mathcal{T}$ and the Kuramoto synchronization. Starting from these observations, scholars developed~\cite{AsllaniCarletti2018,gjata} a control strategy aimed at desynchronize the Kuramoto model based on the Hamiltonian control theory~\cite{vittot,chandre,ciraolo} that ultimately relies on the modification of the Hamiltonian function by the addition of a ``small'' control term, whose goal is to stabilize $\mathcal{T}$ to eventually prevent synchronization in the Kuramoto model.

Let us introduce the reader with few details about the Hamiltonian control theory and referring to, e.g.,~\cite{vittot,chandre,ciraolo} for a more complete study. Let us thus rewrite the Hamilton function $H_0$ as follows, $H_0=H_{int}+V$, where 
\begin{equation*}
	H_{int}=\sum_{i=1}^N\omega_iI_i\text{ and }V=-\frac{K}{N}\sum_{ij=1}^NA_{ij}\sqrt{I_iI_j}(I_i-I_j)\cos(\phi_j-\phi_i)\, ,
\end{equation*}
namely $H_{int}$ is the integrable part and $V$ the ``perturbation'', namely smaller than $H_{int}$. The control term is thus a function $h\thicksim\mathcal{O}(V^2)$ that added to $H_0$ returns a controlled Hamilton function, $H=H_{int}+V+h$, whose dynamics is ``closer'' to the one for $H_{int}$, namely the effects of the perturbation $V$ are reduced. Moreover, the smallness of $h$ makes the control non-invasive. In conclusion, based on the relation between stability of the torus and synchronization of Kuramoto system, by developing a suitable control term we can infer about the dynamical properties of the coupled oscillators.

More precisely, the added control term is defined by 
\begin{equation}\label{eq::whittautControl}
	h:=\sum_{n\geq 1}\frac{\{-\Gamma V\}^n}{(n+1)!}\Gamma\, ,
\end{equation}
where $\{\centerdot\}$ denotes the Poisson bracket and $\Gamma$ is the pseudoinverse operator of $H_{int}$. In~\cite{gjata,AsllaniCarletti2018} authors truncated the previous expression (\ref{eq::whittautControl}) by only retaining the first term. Then rewriting the resulting Hamiltonian equations on the invariant torus, they obtained the following controlled KM 
\begin{equation}
\frac{d{\phi_k}}{dt}=\omega_k+\frac{K}{N}\sum_{j=1}^{N}A_{kj}\sin(\phi_j-\phi_k)+h_k\, ,
\end{equation}
where, for $k=1,\dots,N$,
\begin{equation}\label{eq::h}
	h_k:=-\frac{K^2}{4}\left[R\tilde{R}_k\cos(\psi-\tilde{\psi}_k)-\frac{1}{N}\sum_{l=1}^N\cos(\phi_l-\phi_k)\cos(\tilde{\psi}_l-\phi_l)\tilde{R}_l-\sum_{l=1}^N\frac{\sin(\phi_k-\phi_l)}{\omega_k-\omega_l}\sin(\psi-\phi_l)R\right]\, ,
\end{equation}
$R$ and $\psi$ are defined in (\ref{eq:Order_parameter}) and 
\begin{equation*}
\tilde{R}_ke^{i\tilde{\psi}_k}:=\frac{1}{N}\sum_{j=1}^N\frac{e^{i\phi_j}}{\omega_j-\omega_k}\, .
\end{equation*}
Let us observe that form of the denominators $\omega_j-\omega_k$ is directly related to the pseudoinverse operator $\Gamma$ that imposes a constraint on the proper frequencies.

The latter term, $h_k$, can be simplified~\cite{AsllaniCarletti2018} to obtain a simpler and more effective control term, suitable for numerical analysis and possible implementation. Firstly, the latter two terms in \ref{eq::h} are neglected as they are much smaller than the first one. Then, $\tilde{R}$ and $\tilde{\psi}$ are redefined by only considering information coming from the $M$ controlled nodes, i.e.,
\begin{equation}\label{eq::tildeR}
\tilde{R}_ke^{i\tilde{\psi}_k}:=\frac{1}{M}\sum_{j\in\{i_1,\dots,i_M\}}\frac{e^{i\phi_j}}{\omega_j-\omega_k}\quad \forall k\in\{i_1,\dots,i_M\}
\end{equation}
where $\{i_1,\dots,i_M\}$ is the set of the controller indexes. Finally, to mimic the fact that once an electric signal in injected into a controlled node, an electromagnetic field is also created and it extends beyond the considered node by reaching the nearby ones, hence the control term is assumed to impact neighboring nodes but with an exponential decay as a function of the distance from the controlled node
\begin{equation*}
S^{stim}_k(\phi_{i_1},\dots,\phi_{i_M}):=
\left\{\begin{matrix}
\tilde{h}_k:=-\frac{c K^2}{4}R\tilde{R}_k\cos(\psi-\phi_k)&if~k\in\{i_1,\dots,i_M\}\\
&\\
\sum_{l\in\{i_1,\dots,i_M\}}e^{-2r_{kl}}\tilde{h}_l&\text{ otherwise}
\end{matrix}\right. \, .
\end{equation*}
The parameter $c>0$ characterize the signal strength and $r_{kl}$ is the distance between the nodes $k$ (non controlled node) and $l$ (controlled node), so that the signal received by node $k$ is reduced by an exponential factor depending on the distance. In our numerical analysis we restricted the interactions to the first neighbors, i.e., $r_{kl}=1$.

\section{The functionability}
\label{app:funct}

The functionability has been introduced and studied in~\cite{RosellTarragoDiazGuilera} to have a metric allowing to rank nodes in a networked dynamical system according to their capacity to push the system toward different state, once they are ``perturbed''. More precisely, working with the frustrated Kuramoto model defined on a network, authors have studied how the system state changes, i.e., passes from desynchronized to synchronized one or viceversa, once a phase lag $\alpha$ is added to a single node. Namely they studied the model
\begin{eqnarray*}
\dot{\phi_k}&=&\omega_k+\frac{K}{N}\sum_{j=1}^{N}A_{jk}\sin(\phi_k-\phi_j-\alpha) \\
\dot{\phi_\ell}&=&\omega_\ell+\frac{K}{N}\sum_{j=1}^{N}A_{j\ell}\sin(\phi_\ell-\phi_j) \quad \forall \ell \neq k\, ,
\end{eqnarray*}
where $\alpha\in[0,\pi/2]$ is a fixed parameter; let us notice that the constraint $\alpha < \pi/2$ allows to prove that the system synchronizes. By iteratively considering all the possible $k$, they determined the nodes more or less prone to modify the system behavior. From this intuitive definition, authors obtained a more suitable formula for (a proxy of) the functionability
\begin{equation}
	\label{eq:functionability}
	\mathcal{F}_k(A,\alpha)=\left(\frac{\alpha d_k}{2N^2}\right)^2\sum_{l=1}^{N}\sum_{ij=1}^{N}\left[(L(k,l)L^T(k,l))^{-1}\right]_{ij}\, ,
\end{equation} 
where $d_k$ is the degree of node $k$, $A$ is the network adjacency matrix, $L=\mathrm{diag}(\mathbf{d})-A$ is the Laplace matrix, $\mathbf{d}=(d_1,\dots,d_N)^\top$ the nodes degree vector and $L(k,l)$ is the Laplace matrix whose $k$-th row and $l$-th column have been removed. Therefore, $\mathcal{F}$ is only a function of the network topology and the frustration parameter $\alpha$. In the present work we have considered $\alpha=0.5$ as advised in~\cite{RosellTarragoDiazGuilera}. Interestingly enough authors of~\cite{RosellTarragoDiazGuilera} have shown that more functional nodes have large degree and are peripheral.

Let us conclude by noticing that although functionability is designed to spot the most influential nodes in the global synchronization process for a frustrated Kuramoto model, we decided to apply it to our case as well, even if our goal differs from the one of~\cite{RosellTarragoDiazGuilera}, by relying on previous observations that nodes capable to enhance synchronization seem also good candidate to decrease it. Indeed after all, the functionability could be though to be a property of the underlying network and we found thus interesting to compare it with the pinning scheme in order to investigate its efficiency in this new framework.

\section{The case of the small world networks}
\label{app:ws}
\begin{figure}[ht]
	\centering
	\includegraphics[width=0.5\linewidth]{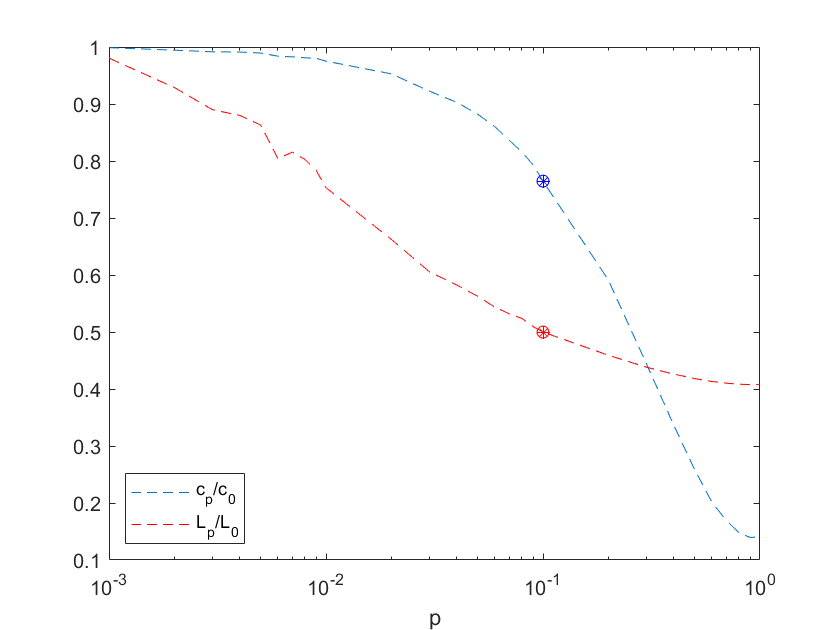}
	\caption{Illustration of the small world effect exhibited by WS networks generated with parameters $N=100$, $\langle k\rangle =10$ and $p_{\textit{WS}}\in[0,1]$. The quantities $L_p$ and $L_0$ represent respectively the average shortest distance when $p_{\textit{WS}}=p$ and $p_{\textit{WS}}=0$; $c_p$ and $c_0$ refer to the average clustering coefficient. The curves were obtained by averaging these quantities over 50 generated networks for each value $p_{\textit{WS}}$. The dots highlight the chosen $p_{\textit{WS}}=0.1$ value.}
	\label{fig::small_world}
\end{figure}

The goal of this section is to apply the control method introduced in the main text to the small-world networks built by using the Watts-Strogatz model (WS)~\cite{WS1998}. Those networks are characterized by high clustering coefficients and small average shortest distance, while they exhibit a quite homogeneous degree distribution. In the following we considered networks formed by $N=100$ nodes, average degree $\langle k\rangle =10$ and rewiring probability $p_{\textit{WS}}=0.1$. Note that this parameter choice enables us to generate graphs that exhibit quite well the small-world effect as state in \cite{WS1998}. This is illustrated in Fig.~\ref{fig::small_world}.

\begin{figure}[ht]
	\includegraphics[width=0.3\linewidth]{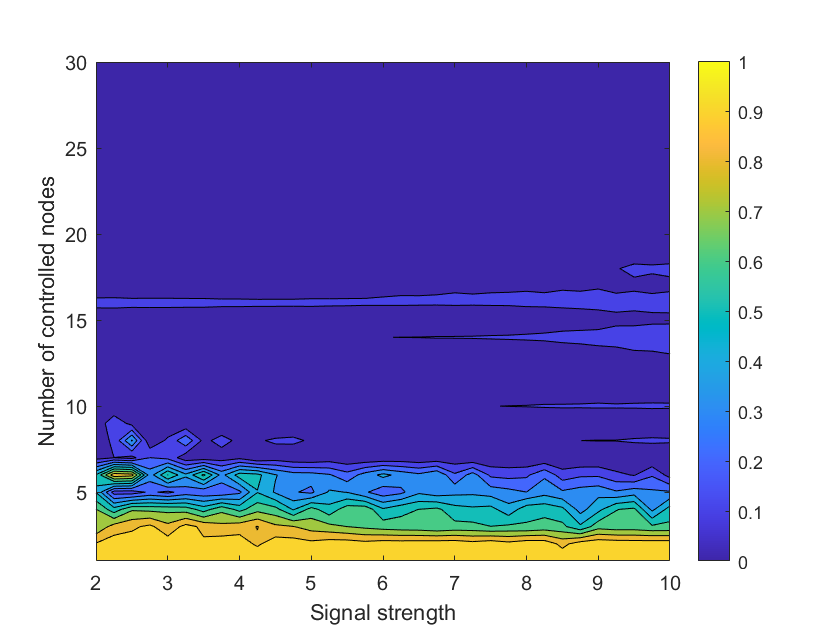}
	\includegraphics[width=0.3\linewidth]{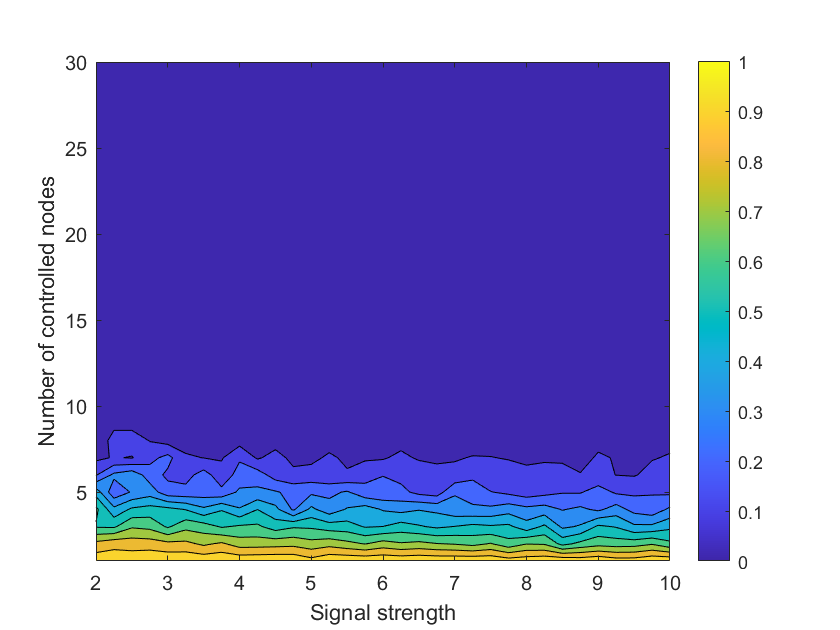}
	\includegraphics[width=0.3\linewidth]{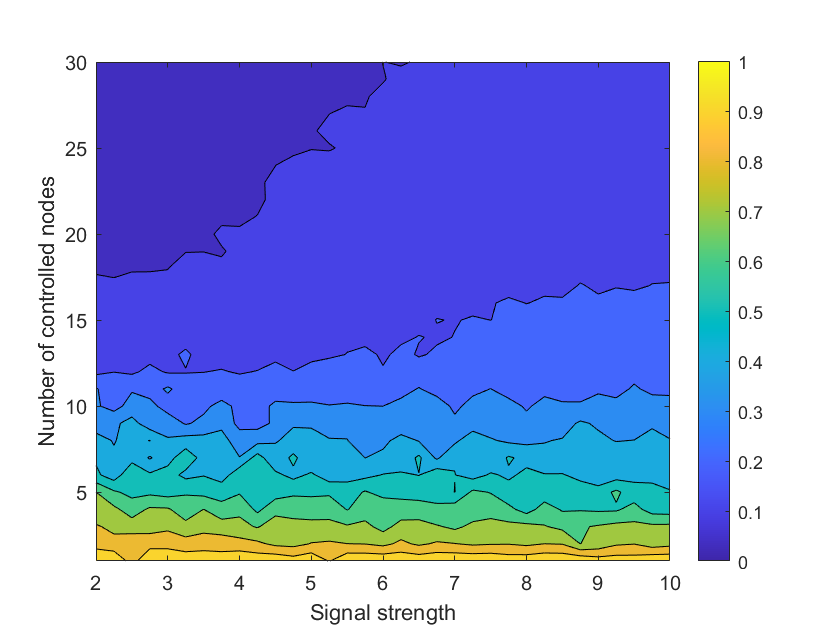}
	\caption{We show the normalized order parameter $\hat{R}$ obtained from numerical simulations of the controlled KM defined on $100$ WS networks composed by $N=100$ nodes, $\langle k\rangle=10$ and probability to rewire $p_{\textit{WS}}=0.1$. Pinned nodes are randomly selected (left panel), according to their degree (middle panel) and the functionability (right panel).}
	\label{Fig:sim_WS}
\end{figure}

\begin{figure}[ht]
	\centering
	\includegraphics[width=0.49\linewidth]{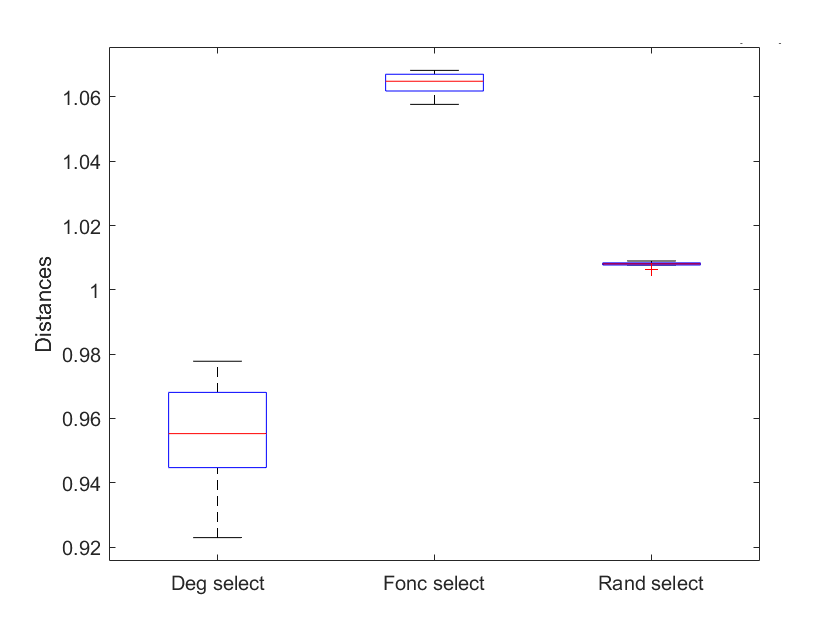}
	\includegraphics[width=0.49\linewidth]{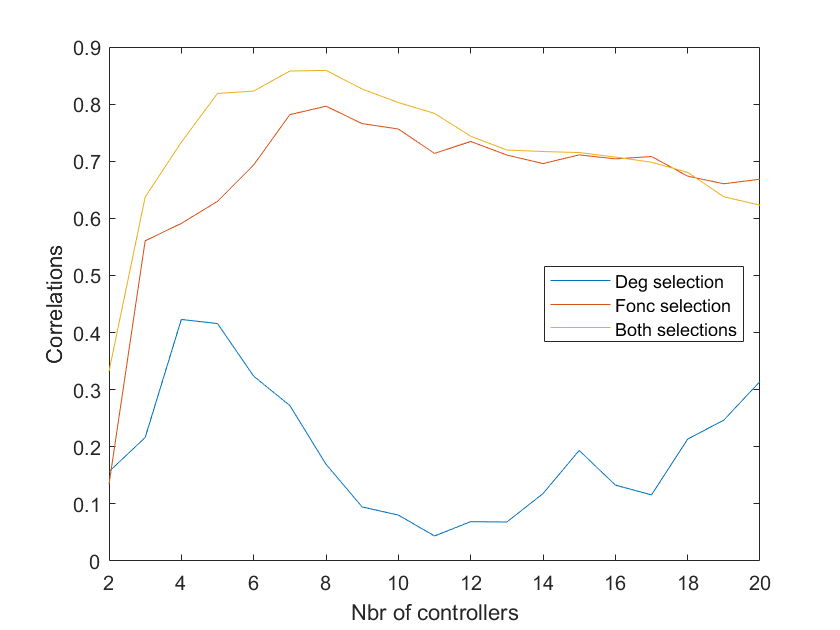}
	\caption{Left panel: distribution of the average normalized distances between controllers and the remaining nodes in the WS networks. We can observe that degree-based selection of controllers tends to decrease those distances, while the opposite is observed with the functionability selection strategy. Right panel : we report the correlation between $\hat{R}$ and above mentioned distance in function of the number of controllers for the degree-based and functionability selection methods. The correlation is quite low in the case of degree-based selection (blue curve) and high with functionability (red curve). In general, it is strongly positively correlated (we can see that by considering all simulations together regardless of the selection methods, yellow curve), in particular when the number of controllers is small.}
	 \label{Fig:distances}
\end{figure}

In Fig.~\ref{Fig:sim_WS} we report the normalized order parameter $\hat{R}$ computed by using $100$ WS networks with the parameters given above, as a function of the number of controllers and the signal strength. One can observe that functionability-based selection provides less good results than both the others; indeed, the dark blue area on the right panel, corresponding to couples $(M,\gamma)$ giving rise to $\hat{R}\leq 0.15$, is smaller than in the other two cases. Those results, at first glance surprising, can be explained by the fact that nodes in WS networks do exhibit very similar degree and moreover they are close each other, hence there are not node with very large or very small betweenness. In Fig.~\ref{Fig:distances} we report the distribution of distances between unpinned and pinned nodes, selected according to the three selection strategies, let us note that such distances are normalized by dividing by the average distance computed on the whole network. We can observe that in the case of functionability those distances are slightly larger than in the remaining cases, moreover, there is a highly positive correlation between those distances and the normalized order parameter $\hat{R}$ (see Fig.~\ref{Fig:distances}). It seems thus that the functionability tends to select node that are quite clustered together in small-world networks and not well connected to the others, this feature could be one of the reasons that make the control less efficient. 

As with star network, we can so observe a similarity with the results of~\cite{LuLiRong2009}. In the latter, authors stated that there is a correlation between the minimal nonzero eigenvalue of the reduced Laplacian matrix $\tilde{L}$, obtained by removing the lines and columns corresponding to the controllers indexes, and the controllers subset efficiency. By simulating dynamical systems of networked linearly coupled oscillators, authors assessed that the larger this eigenvalue, the better the system synchronizes. Despite the fact that our system exhibits non-linear coupling function and that our goal is, in opposition, to get the system desynchronize, an interesting parallel could be drawn.

\begin{figure}[ht]
	\centering
	\includegraphics[width=0.4\linewidth]{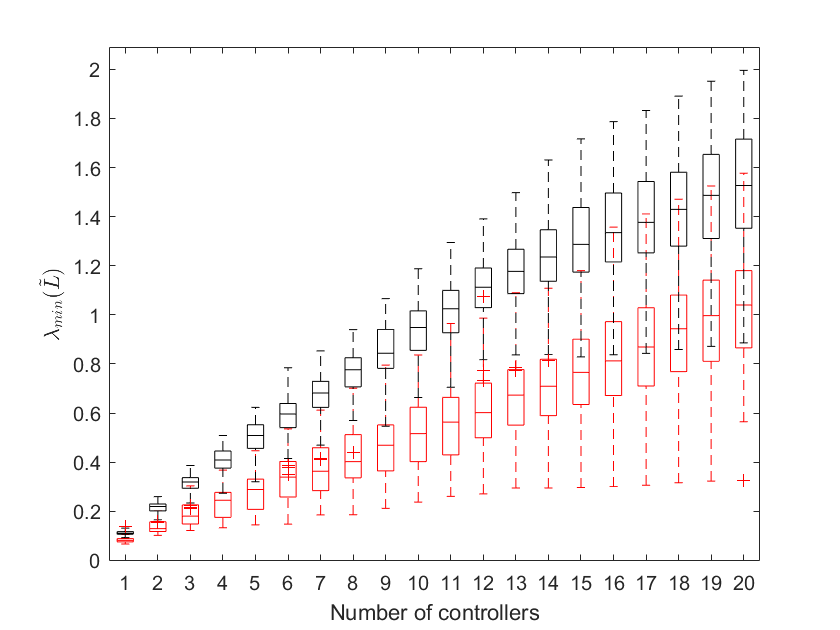}
	\caption{Boxplots of values of $\lambda_{min}(\tilde{L})$ among WS networks in function of the number of controllers selected according to their degree (black boxes) and functionability (red boxes).}
	\label{Fig::wseig}
\end{figure}

We computed minimal eigenvalues of the corresponding matrices $\tilde{L}$ and show them in Fig.~\ref{Fig::wseig}. From the figure it is clear that $\lambda_{min}(\tilde{L})$ are larger when pinned nodes are selected according to their degree than functionability, and this strengthen the parallel drawn with the linear synchronization case.

\section{Betweenness centrality}
\label{app:betweenness}
\begin{figure}[ht]
	\centering
	\includegraphics[width=0.4\linewidth]{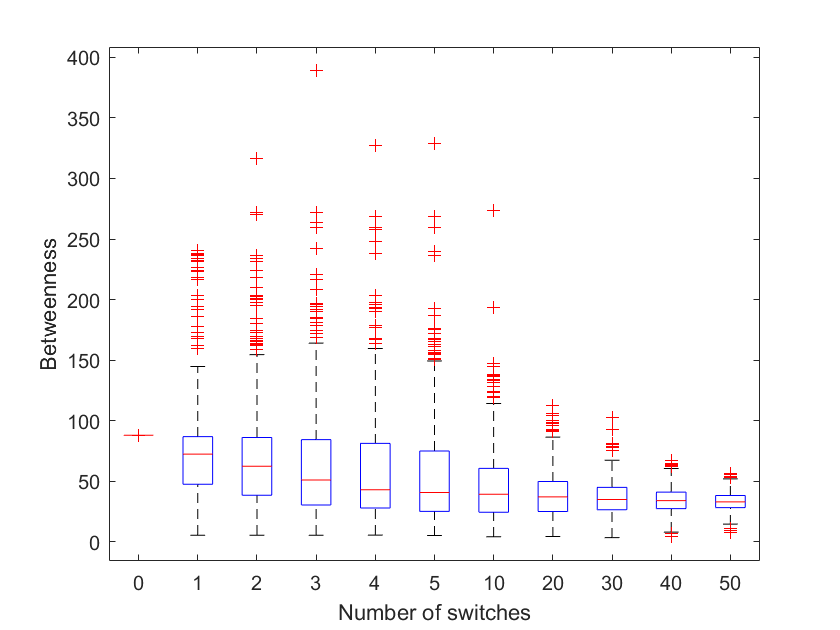}
	\caption{We report here the betweenness centrality distribution among networks of $N=50$ nodes in function of the number of edge switches (starting from a regular ring-shaped lattice in which all nodes have 6 neighbors, each switch exchanges extremities of two edges that do not share any incident node). Note that the regular degree distribution is unchanged through rewiring.}
	\label{fig::between_disp}
\end{figure}

The aim of this section is to present and discuss the numerical results obtained by using the betweenness centrality to select pinned nodes. As we stated in Appendix~\ref{app:funct}, the functionability ranks higher, nodes with large degree or low betweenness, we are thus interested in studying the impact of the latter on the synchronization. However to disentangle the possible contributions from these two factors we created networks starting from regular lattices with periodic boundary conditions, i.e., whose nodes have all the same degree, and then we reshuffle the connection by keeping constant the degrees sequence (criss-cross reshuffling). In particular we used regular lattices with degree $k_i=6$ composed by $N=50$ nodes and we performed $5$ links switches. This very last choice is motivated by an analysis of the range of values for the betweenness. As we can see in Fig. \ref{fig::between_disp}, only a small number of switches are required to have a quite large variation in the betweenness values. We fix here the choice at 5 to be close of the maximal variation while having still a large diversity of networks shapes.

\begin{figure}[ht]
	\centering
	\includegraphics[width=0.45\linewidth]{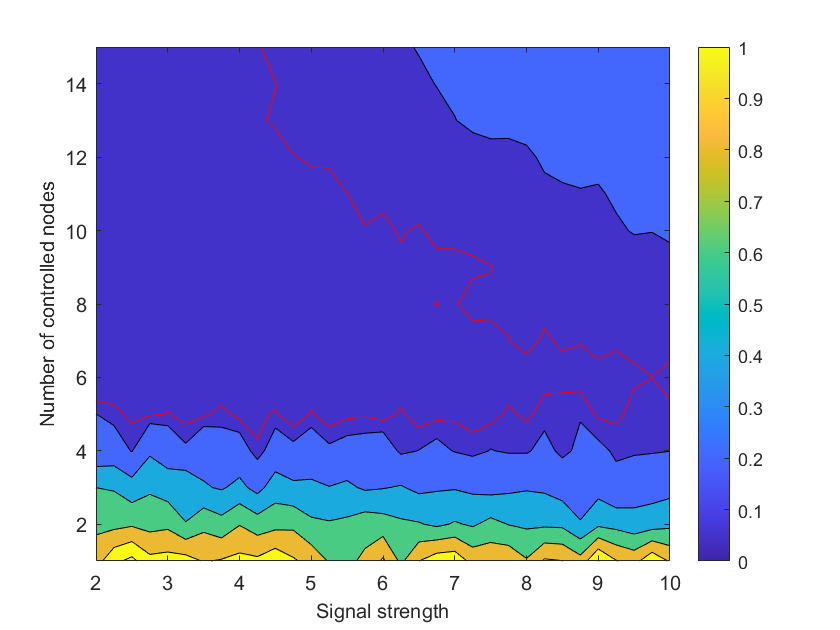}
	\includegraphics[width=0.45\linewidth]{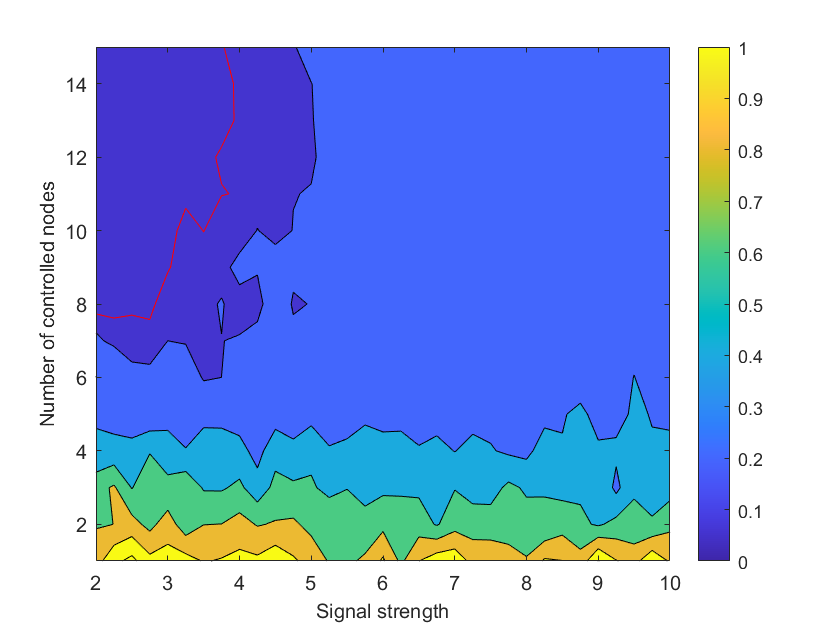}
	\caption{We report the normalized order parameter $\hat{R}$ obtained from numerical simulations of the controlled KM applied on $100$ networks built using the procedure described in the text. On the left panel, the controllers were selected by using large betweenness, the opposite was done for the right panel.}
	\label{fig::between}
\end{figure}

To check the impact of the betweenness centrality we compared two opposite selection strategies, one that favors the selection of nodes with high betweenness and a second that, on the opposite, select preferentially nodes with low betweenness. The results reported in Fig.~\ref{fig::between} show that the synchronization level is higher for the second strategy. They therefore support the claim that small betweenness is not a relevant factor and this fact can partially explain the bad results of functionability selection in the case of graphs where degree is homogeneous. 

\section{Some additional numerical results for the core-periphery case.}
\label{sec:figureapp}

\begin{figure}[!t]
	\centering
	\includegraphics[width=0.3\linewidth]{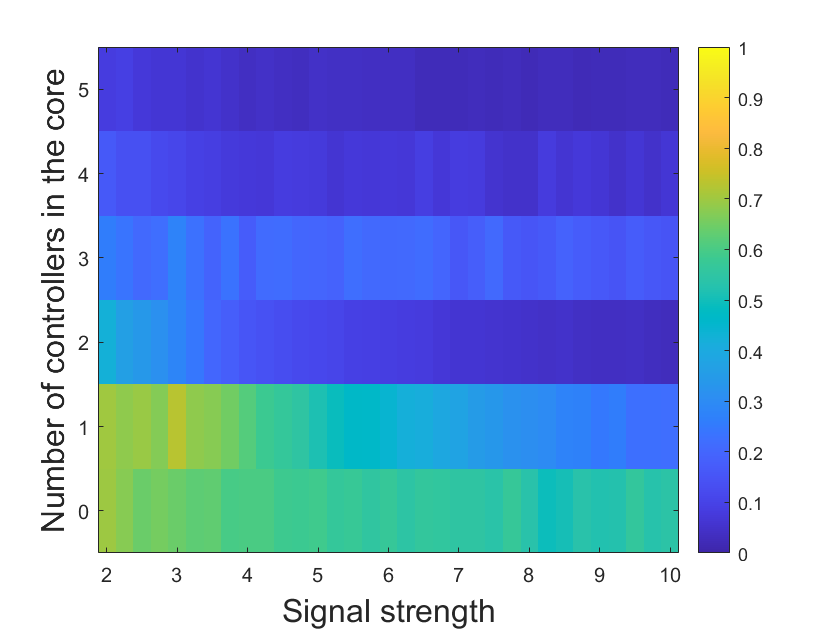}
	\caption{We show the normalized order parameter $\hat{R}$ obtained from numerical simulations of the controlled KM applied on $100$ core-periphery networks with $(n,M,p)=(10,5,0.7)$. Controllers in the core and in the periphery were selected uniformly at random.}
	\label{Fig::cp10App} 
\end{figure}

\begin{figure}
	\centering
	\includegraphics[width=0.3\linewidth]{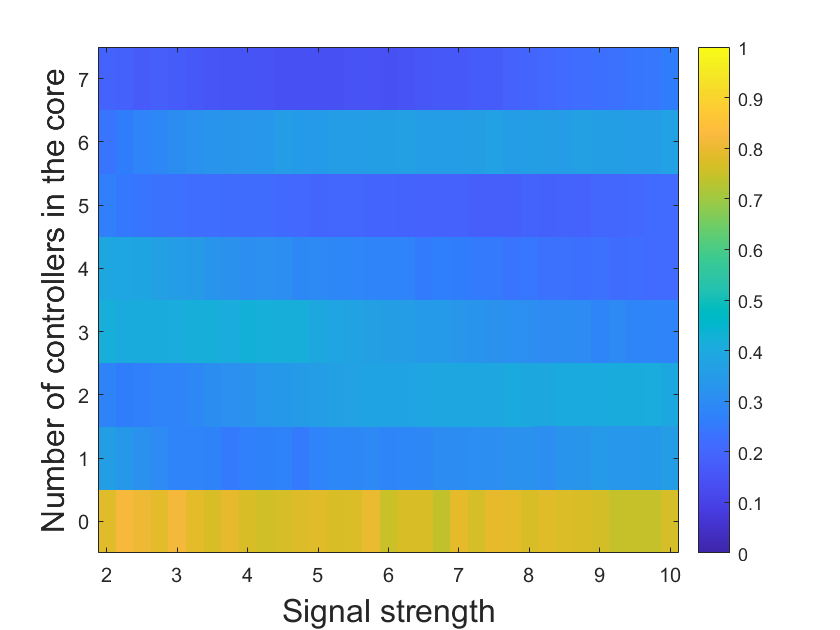}
	\includegraphics[width=0.3\linewidth]{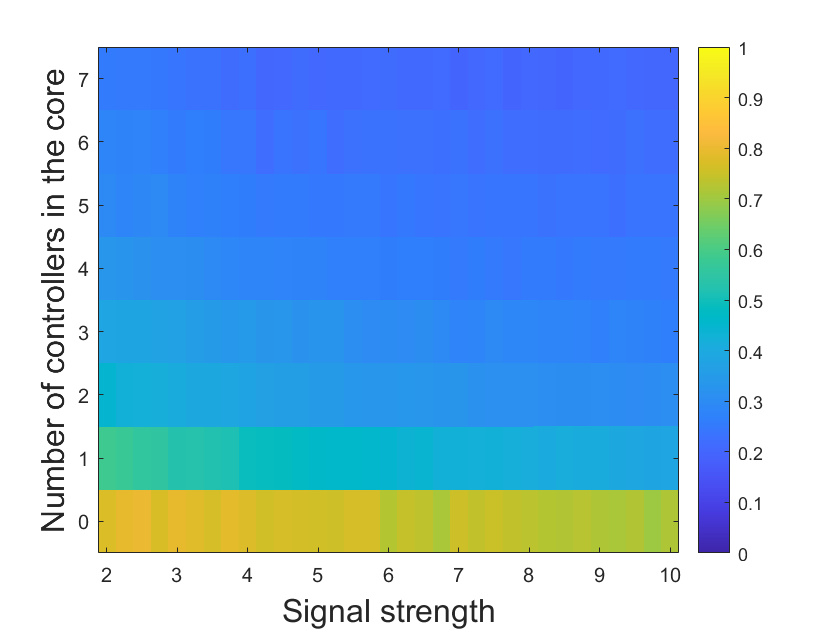}
	\includegraphics[width=0.3\linewidth]{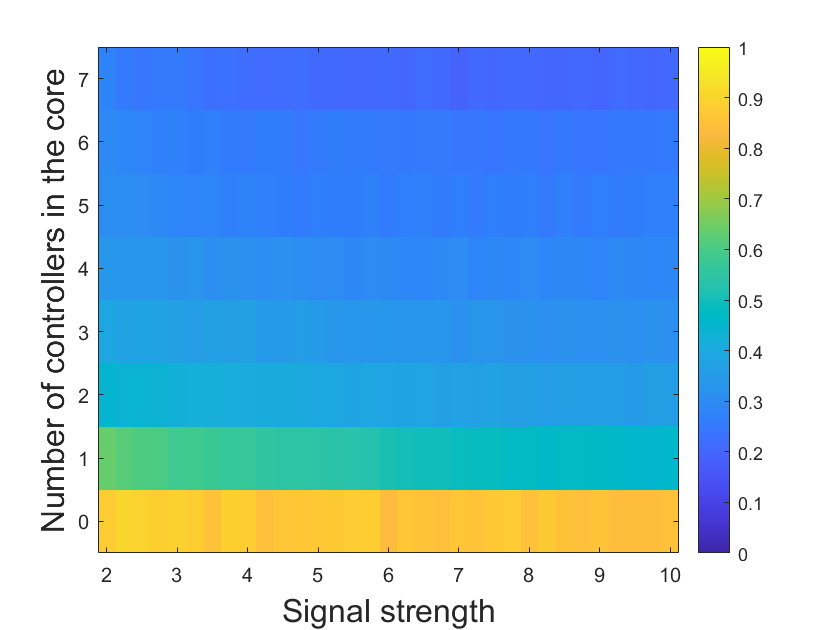}
	\caption{We show the normalized order parameter $\hat{R}$ obtained from numerical simulations of the controlled KM applied on $100$ core-periphery networks with $(n,M,p)=(20,7,0.7)$. Controllers were selected according to three selection methods, uniformly at random (left panel), degree-based selection (middle panel) and functionability-based selection (right panel).}
	\label{Fig::cp20} 
\end{figure}

The aim of this section is to present some additional results in the core-periphery that bring more details and strengthen the conclusions we provided in the main part of the paper. In Fig.~\ref{Fig::cp10App} we show the normalized order parameter for the core-periphery network with $(n,M,p)=(10,5,0.7)$ where the controllers, although we impose that some belong to the core, are selected uniformly at random in the core and periphery. In Fig.~\ref{Fig::cp20} we present the case $(n,M,p)=(20,7,0.7)$. We can clearly observe that both case exhibit a similar behavior and thus the same kind of conclusions can be drawn.

\begin{figure}
	\centering
	\includegraphics[width=0.75\linewidth]{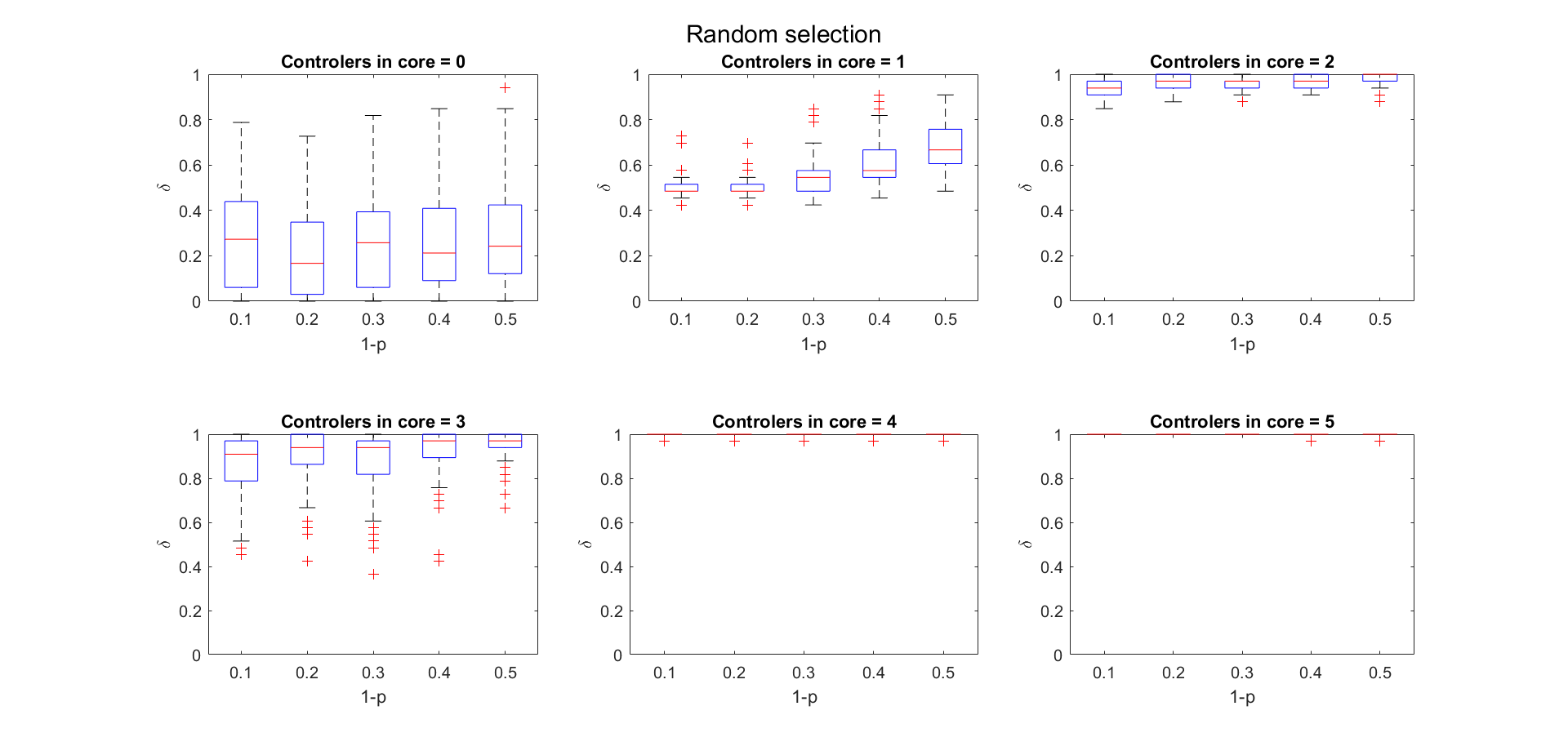}
	\includegraphics[width=0.75\linewidth]{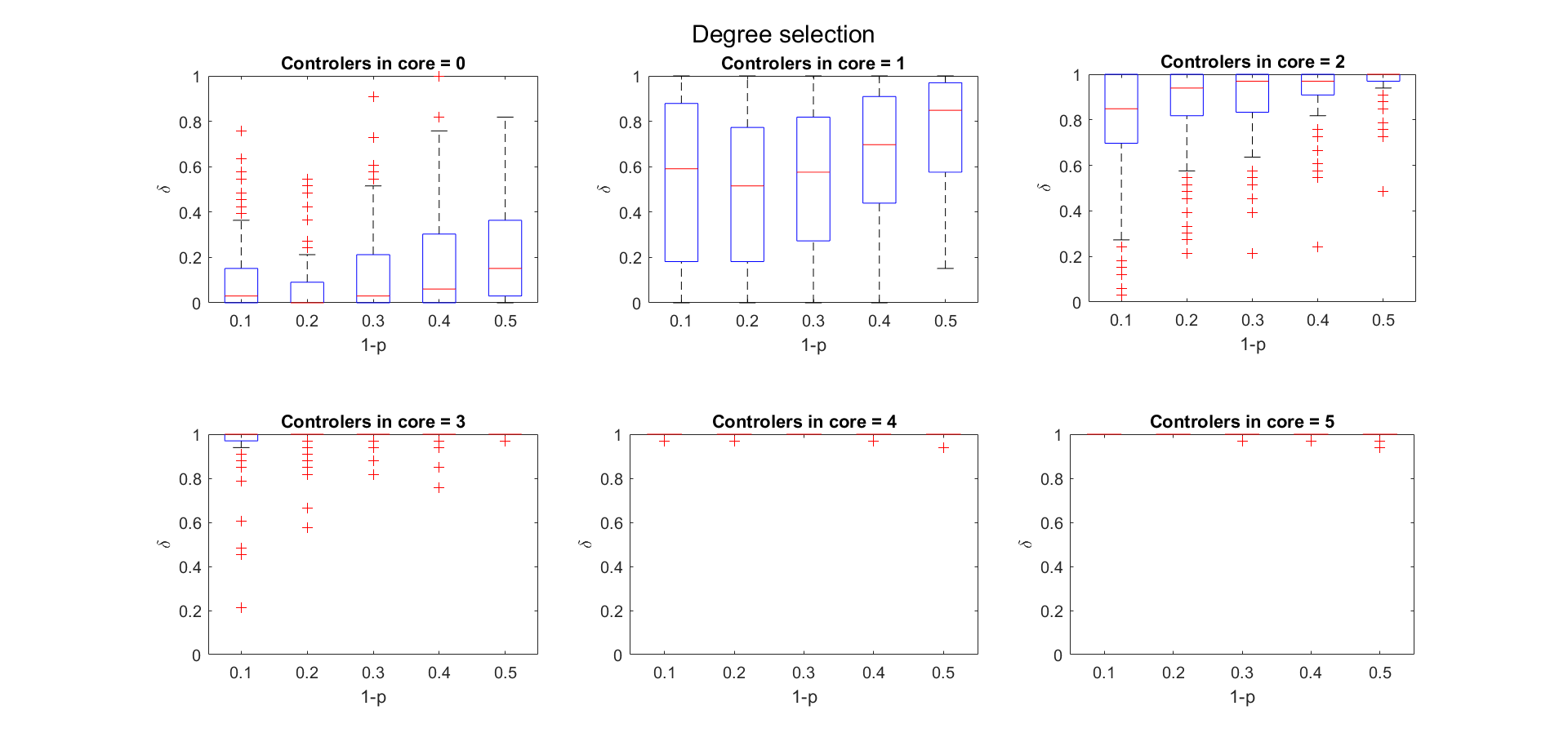}
	\includegraphics[width=0.75\linewidth]{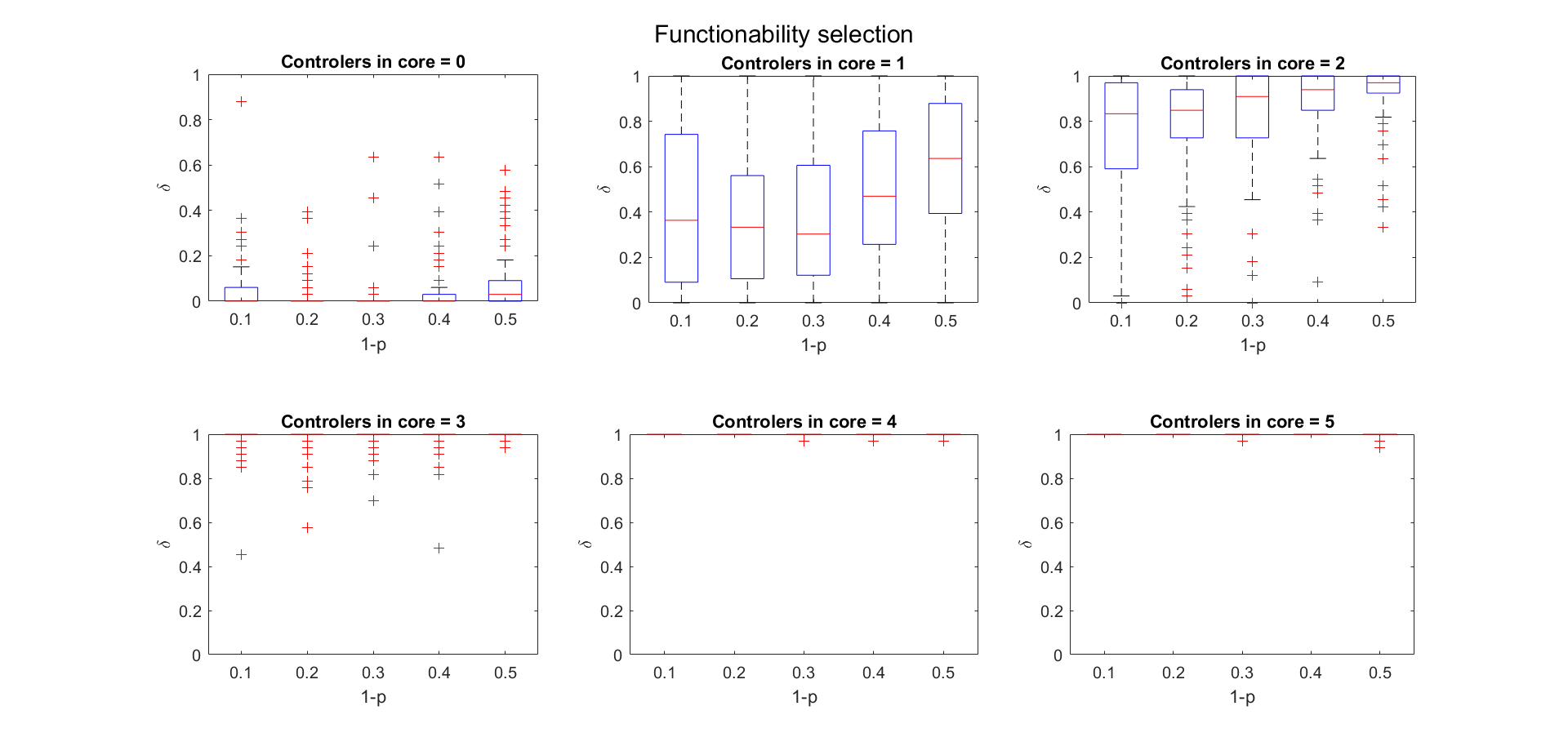}
	\caption{We show the distribution of the size, $\delta$, of parameters region resulting into a well desynchronized state, i.e., associated to $\hat{R}\leq 0.15$, for different core-periphery networks made by $N=100$ nodes, generated for several values of $p$. Panels in the top two rows: pinned nodes are randomly selected; third and fourth middle panels, nodes are selected according to their degree; fifth and sixth bottom panels: selection is performed according to functionability. For each selection strategy we consider a different assignment of controllers in the core and in the the periphery, as reported in the panel titles.}
	\label{Fig::app::cp_p} 
\end{figure}

Finally, in Fig.~\ref{Fig::app::cp_p} we report a more detailed version of results already shown in Fig.~\ref{Fig::cp_p}, by using in particular the three selection strategies, random, node-degree and functionability, and also a finer grid for the number of controlled nodes in the core, i.e., $M=0,\dots,5$. We can observe that random selection performs worst than the other two strategies except for $M=4$ and $M=5$ controlled nodes in the core, and thus only $1$ and $0$ in the periphery, however in these last two cases the three methods all perform very well by returning very large values of $\delta$. The degree based selection is more efficient than the functionability one for $M=0,1,2$, the two methods are almost equivalent for $M=3,4,5$.

\section{Simulations on reshuffled real networks}
\label{sec:reshuffled}

Let us conclude by considering some supplementary simulations that we performed using the three real networks already presented in the main text but applying a degree distribution preserving process of links reshuffling; the aim of this analysis is to check if the obtained desynchronization results depend on the degree or on other ``hidden" topological features of the network. Hence for a given network, we chose two links uniformly at random (they must be incident to four different nodes) and we rewire them in order to exchange one of their extremities, this process has been called ``criss-cross" because the links swap recall a cross. 

\begin{figure}[!h]
	\includegraphics[width=0.33\linewidth]{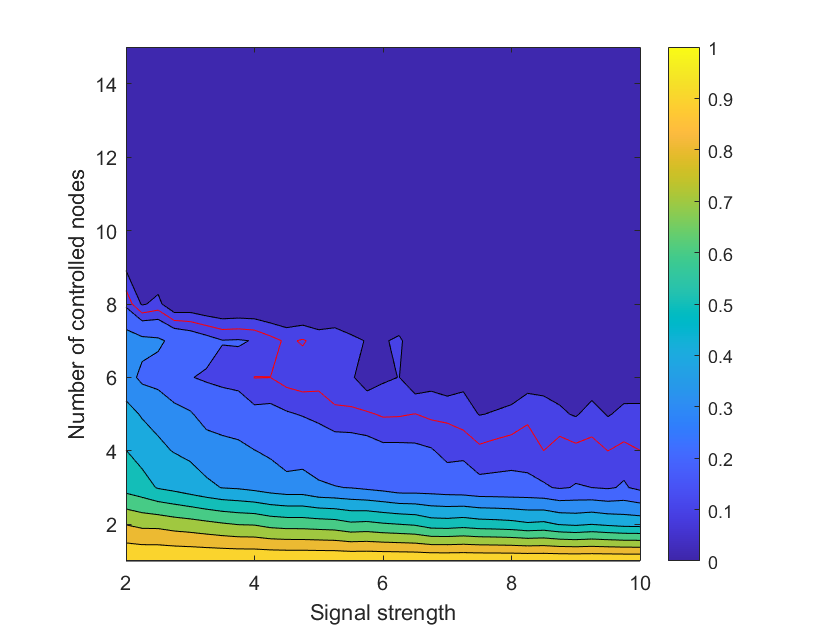}\hfill
	\includegraphics[width=0.33\linewidth]{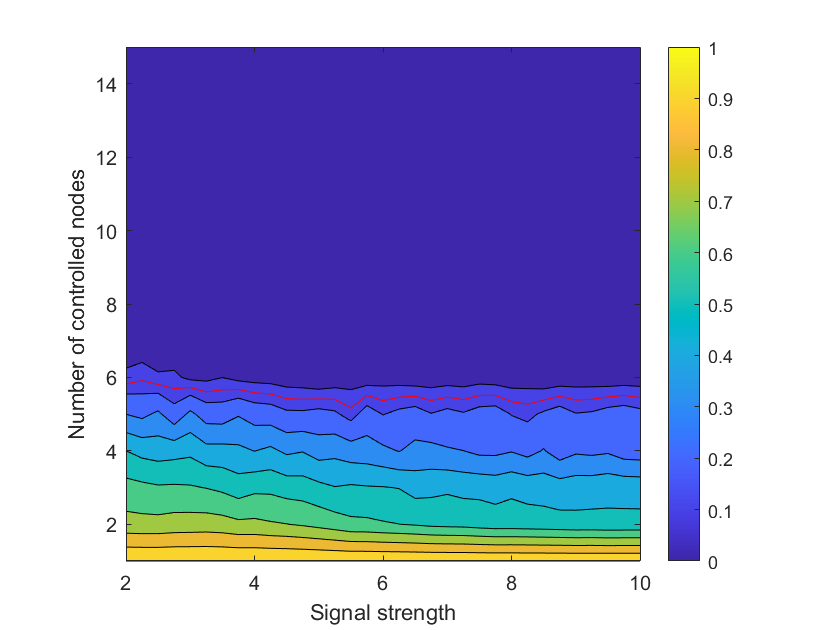}\hfill
	\includegraphics[width=0.33\linewidth]{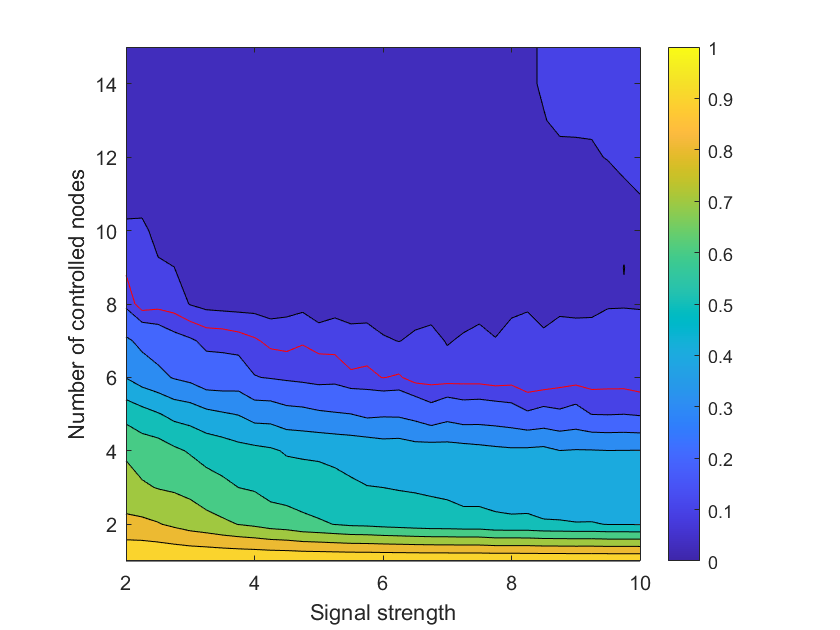}\\
	\includegraphics[width=0.33\linewidth]{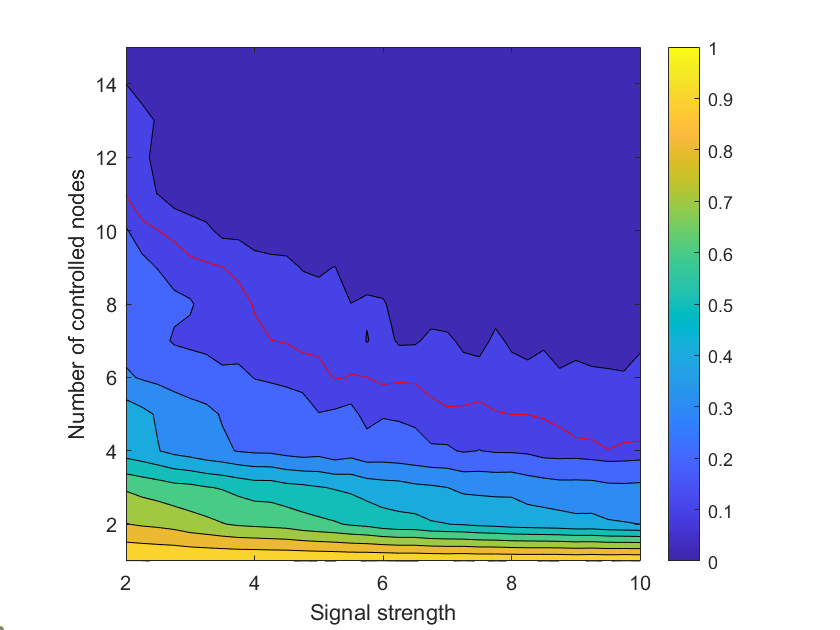}\hfill
	\includegraphics[width=0.33\linewidth]{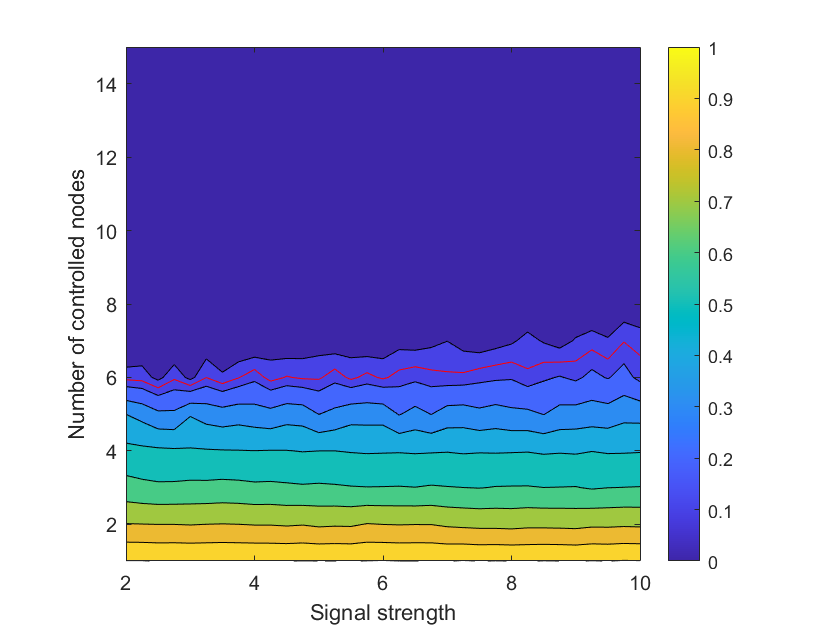}\hfill
	\includegraphics[width=0.33\linewidth]{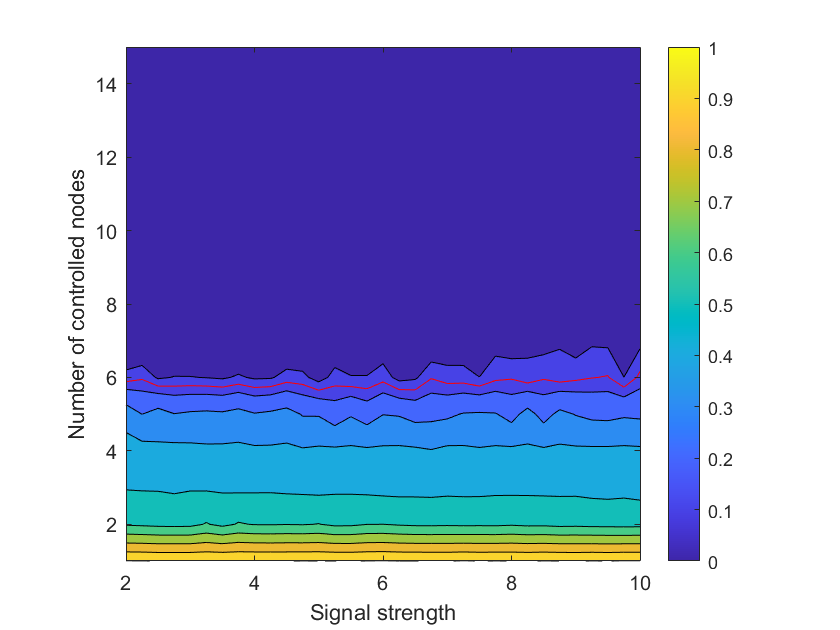}\\
	\includegraphics[width=0.33\linewidth]{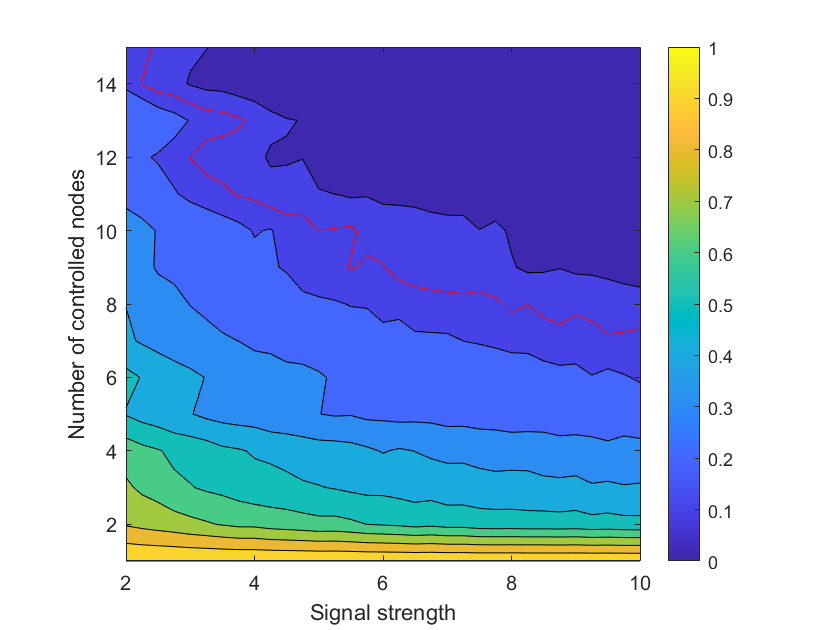}\hfill
	\includegraphics[width=0.33\linewidth]{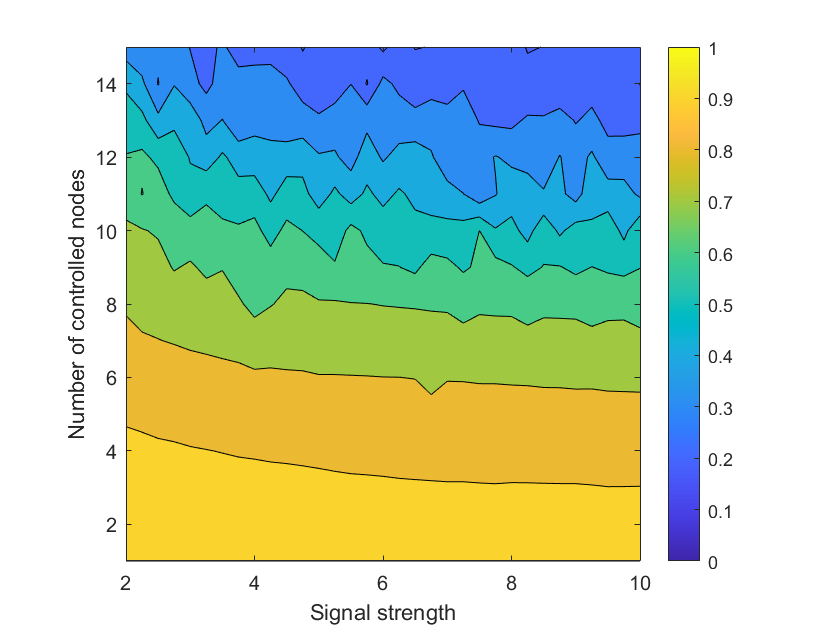}\hfill
	\includegraphics[width=0.33\linewidth]{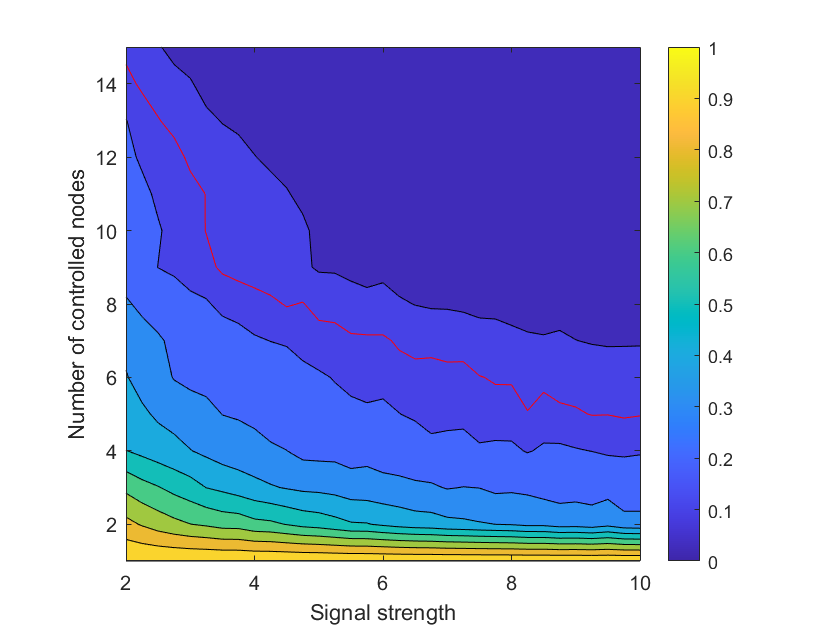}\\
	\caption{Normalized order parameter $\hat{R}$ obtained from the cat connectome network in which $L/2$ links rewiring preserving the degree distribution have been made. Each row corresponds to one connectome (First: cat; second: macaque; third: mouse) and each column to one selection strategy (First: uniformly at random; second: degree; third: functionability).\label{fig::real_networks_reshuffle}}
\end{figure}

Starting from the connectome networks presented above, we applied the rewire process $\frac{L}{2}$ times, being $L$ the number of links in the network, so that each link has the chance to be rewired. Let us stress again that the network topology is modified but not its degree distributions. By using those rewired networks, we computed the average normalized order parameter, $\hat{R}$, as explained in the main text. The results obtained with those modified networks (see Fig.~\ref{fig::real_networks_reshuffle}) are very close to the ones obtained for the original networks (see Fig.~\ref{fig::real_networks}). This observation can be put on a quantitative ground by looking at Table~\ref{table:sim_realAPP}, where we present the fraction, $\delta$, of the parameters values (number of controllers and control signal strength) for which $\hat{R}\leq 0.15$, for both the original and the reshuffled networks. It can be observed that the differences are of the order of $10^{-2}$. 
\begin{table}[ht]
	\centering
	\begin{tabular}{|c|c|c|c|c|}
		\hline
		&&\textbf{Random selection}&\textbf{Degree selection}&\textbf{Functionability selection}\\
		\hline
   		\multirow{3}*{\textbf{Cat}}
   		&\textbf{Original}&0.6788    &0.6808    &0.5556\\
		&\textbf{Rewired}&0.6606    &0.6667    &0.6121\\
		\hline
		\multirow{3}*{\textbf{Macaque}}
		&\textbf{Original}&0.6101    &0.6646    &0.6667\\
		&\textbf{Rewired}&0.5980    &0.6263    &0.6626\\
		\hline
		\multirow{3}*{\textbf{Mouse}}
		&\textbf{Original}&0.3636    &0.0162    &0.5152\\
		&\textbf{Rewired}&0.3778    &     0    &0.5232\\
		\hline
	\end{tabular}
	\caption{
	We report the fraction $\delta$ of parameters, number of controllers and control signal strength, for which $\hat{R}<0.15$ by running the controlled KM on the three real connectome networks (namely to Fig.~\ref{fig::real_networks}) and their rewired versions (namely to the Fig.~\ref{fig::real_networks_reshuffle}).}
	\label{table:sim_realAPP}
\end{table}

The last observation is interesting and suggests that if two networks have the same degree distribution, then they exhibit a similar desynchronization behavior. This is however only a preliminary observation and some future work is needed to tackle the question. Was it true, than it would be enough to know the degree distribution of the network one wants to control, then generate an artificial network with the same degree distribution and test pinning control schemes on this network to learn how to control the original one. 


\begin{thebibliography}{99}

\bibitem{pikovsky2003} Pikovsky, A., Rosenblum, M. and Kurths, J. , \textit{Synchronization: A Universal Concept in Nonlinear Sciences}, \textbf{12}, Cambridge University Press, Cambridge, (2003).

\bibitem{goldbeter_berridge_1996} Goldbeter, Albert and Berridge, M. J., \textit{Biochemical Oscillations and Cellular Rhythms: The Molecular Bases of Periodic and Chaotic Behaviour}, Cambridge University Press, Cambridge (1996)

\bibitem{kuramoto1975} Kuramoto, Y. , \textit{Self-entrainment of a population of coupled non-linear oscillators} in \textit{Lecture Notes in Physics, International Symposium on Mathematical Problems in Theoretical Physics} (Araki, H.) \textbf{39}, 420--422 (1975).

\bibitem{kuramoto} Kuramoto, Y. , \textit{Chemical Oscillations, Waves, and Turbulence}, Springer-Verlag (1984).

\bibitem{gjata} Gjata, O., Asllani, M., Barletti, L. and Carletti, T. , \textit{Using Hamiltonian control to desynchronize Kuramoto oscillators}, \textit{Phys. Rev. E} \textbf{95}(2), 022209 (2017).

\bibitem{AsllaniCarletti2018} Asllani, M., Expert, P.  and Carletti, T.,  PLoS Comput Biol, $\mathbf{14}$ (7), (2018), pp. e1006296

\bibitem{vittot} Vittot, M. Perturbation theory and control in classical or quantum mechanics by an inversion formula. \textit{J. Phys. A: Math. Gen.} \textbf{37}, 6337 (2004).

\bibitem{chandre} Chandre, C. \textit{et al.} Channeling chaos by building barriers. \textit{Phys. Rev. Lett.} \textbf{94}, 074101 (2005).

\bibitem{ciraolo} Ciraolo, G. \textit{et al.} Control of Hamiltonian chaos as a possible tool to control anomalous transport in fusion plasmas. \textit{Phys. Rev. E} \textbf{69}, 056213 (2004).

\bibitem{carletti} Boreux, J., Carletti, T., Skokos, C., Papaphilippou, Y. and Vittot, M. Efficient control of accelerator maps. \textit{Int. J. Bifurcation Chaos} \textbf{22}, 1250219 (2012).

\bibitem{pyragas} Pyragas, K., Popovych, O. and Tass, P. Controlling synchrony in oscillatory networks with a separate stimulation-registration setup. \textit{EPL}, \textbf{80}, 40002 (2007). 

\bibitem{ER} Erd\H{o}s, P. and R\'enyi, A. On Random Graphs. I, Publicationes Mathematicae, $\mathbf{6}$ (3-4), 290 (2022)

\bibitem{Bollobas} Bollob\'as, B. Random Graphs (2nd ed.) Cambridge University Press (2001). 

\bibitem{NW} Newman, M. E. J. and Watts, D. J. Renormalization group analysis of the small-world network model, Physics Letters A, $\mathbf{263}$, 341, (1999).

\bibitem{Parkinson} Sveinbjornsdottir, S. The clinical symptoms of Parkinson's disease. \textit{J. Neurochem.} \textbf{139}, 318 (2016).

\bibitem{Parkes:1974wr} J. D. Parkes, R. C. Baxter, C. D. Marsden, and J. E. Rees, J. Neurol. Neurosurg. Psychiatry \textbf{37}, 422 (1974).

\bibitem{Holloway:2005wh} R. G. Holloway, Archives of Neurology \textbf{62}, 430 (2005).

\bibitem{Deckers:2000tz} C. Deckers, S. J. Czuczwar, Y. A. Hekster, A. Keyser, H. Kubova, H. Meinardi, P. N. Patsalos, W. O. Renier, and C. M. Van Rijn, Epilepsia \textbf{41}, 1364 (2000).

\bibitem{neurostim} Hallet, M. Transcranial magnetic stimulation and the human brain. \textit{Nature} \textbf{406}, 147 (2000).

\bibitem{tass} Tass, P. Desynchronization by means of a Coordinated Reset of neural sub-populations: A novel technique for demand-controlled Deep Brain Stimulation. \textit{Prog. Theor. Phys. Suppl.} \textbf{150}, 281 (2003).

\bibitem{neurocontrol} Popovych, O. and Tass, P. Control of abnormal synchronization in neurological disorders. \textit{Front. Neurol.} \textbf{5}, 268/fneur.2014.00268 (2014).

\bibitem{DBS} Kringelbach, M., Jenkinson, N., Owen, S. and Aziz, T. Translational principles of deep brain stimulation. \textit{Nat. Rev. Neurosci.} \textbf{8}, 623 (2007).

\bibitem{Bronstein:2011hh} Bronstein J. M. \textit{et. al.}, \textit{Arch. Neurol.} \textbf{68}, 1 (2011).

\bibitem{Theodore:2004gv} Theodore W. H. and Fisher R. S., \textit{Lancet Neurol.} \textbf{3}, 111 (2004).

\bibitem{RosellTarragoDiazGuilera} Rosell-Tarrag\'o, G. And D\'{\i}az-Guilera, A., Functionability in complex networks: Leading nodes for the transition from structural to functional networks through remote asynchronization, Chaos \textbf{30}, 013105 (2020)

\bibitem{Witthaut} Witthaut D. and Timme M., Kuramoto dynamics in Hamiltonian systems. Phys. Rev. E 90, 032917 (2014)

\bibitem{WangChen2002} Wang X. and Chen G., Pinning control of scale-free dynamical networks, Physica A: Statistical Mechanics and its Applications, 310, p 521-531, (2002)

\bibitem{LuLiRong2009} Lu W., Li X. and Rong Z., Global stabilization of complex networks with digraph topologies via a local pinning algorithm, Elsevier: Automatica 46 (2010) 116-121

\bibitem{Liu2018} H. Liu, X. Xu, J. -A. Lu, G. Chen and Z. Zeng, "Optimizing Pinning Control of Complex Dynamical Networks Based on Spectral Properties of Grounded Laplacian Matrices,"  IEEE Transactions on Systems, Man, and Cybernetics: Systems, vol. 51, no. 2, pp. 786-796, Feb. 2021, doi: 10.1109/TSMC.2018.2882620.

\bibitem{Borgatti} Borgatti S. \& Ervertt M., Models of core/periphery structure, Elsevier: Social Networks 21, p 375-395, 1999.

\bibitem{ZhouMondragon} Zhou S. \& Mondragon R., The rich-club phenomenon in th internet topology, IEEE Communications Letters, vol. 8, no. 3, pp.180-182, March 2004

\bibitem{Colizza} Colizza V., Flammini A., Serrano M. A. \& Vespignani A., Detecting rich club ordering in complex networks, Nature Physics 2, 110-115 (2006)

\bibitem{vdheuvel} van den Heuvel M. \& Sporns O.  Rich-Club Organization of the Human Connectome, Journal of Neuroscience 2 November 2011, 31 (44) 15775-15786

\bibitem{LNR2017} V. Latora, V. Nicosia, G. Russo, {\em Complex Networks: Principles, Methods and Applications}, Cambridge University Press, (2017)

\bibitem{WS1998} Watts, Duncan J. and Strogatz, Steven H., {\em Collective dynamics of ‘small-world’ networks}, Nature, $\mathbf{393}$, (1998), pp. 440.

\bibitem{Barthelemy} Barthelemy M., Spatial Networks, Physics Reports, \textbf{499}, 1-101 (2011) 

\bibitem{networkswebsite} CONNECTOMES: Neurodata's Graph Database, Neurodata, \url{https://neurodata.io/project/connectomes/}, last access on the 14th of February 2023.

\end{thebibliography}
\end{document}